\documentclass[oneside]{article}
\usepackage[a4paper, total={6in, 8in}]{geometry}
\usepackage{bm}
\usepackage{amsmath}
\usepackage{hyperref}
\usepackage{multicol}
\usepackage{amsfonts}
\usepackage{graphicx}
\usepackage{amsthm}
\usepackage{cleveref}
\usepackage{subcaption}
\usepackage{mathrsfs}
\usepackage{authblk}
\usepackage{float}
\usepackage{xcolor}
\usepackage{amssymb} 
\usepackage{amsmath} 
\usepackage{xcolor}
\usepackage{bm} 
\usepackage{newtxtext}
\usepackage{newtxmath}
\usepackage{graphicx}
\usepackage{subcaption}
\usepackage{sidecap} 

\usepackage[LGR,T1]{fontenc}

\usepackage{amsmath, amsthm, amssymb, verbatim, comment, tikz, pgfplots, mathtools, mathrsfs, ifthen, import, float, subcaption, enumitem, todonotes, bbm}
\tikzset{>=latex} 
\usepackage{tikz-network}
\usetikzlibrary{decorations.pathreplacing}
\usepackage{pgfplots}

\usepackage{amsfonts}
\usepackage[ruled,vlined,linesnumbered]{algorithm2e}
\usepackage{fullpage}
\usepackage{xifthen}
\usepackage{enumerate}
\usepackage{titlesec}
\newcounter{reviewer}
\newcounter{point}[reviewer]
\newtheorem{theorem}{Theorem}[section]
\newtheorem{remark}[theorem]{Remark}
\newtheorem{proposition}[theorem]{Proposition}

\newtheorem{corollary}[theorem]{Corollary}

\begin{document}
\title{\bf{A Least-Squares-Based Neural Network (LS-Net) for Solving Linear Parametric PDEs}}

\author{Shima Baharlouei\textsuperscript{1}\thanks{Corresponding author, email:shima.baharlouei@ehu.eus}\hspace{0.2cm},~Jamie M. Taylor\textsuperscript{2},~Carlos Uriarte\textsuperscript{1},~David Pardo\textsuperscript{1, 3, 4}\\%

\textsuperscript{1}{\small{University of the Basque Country (UPV/EHU), Leioa, Spain.}}\\%
\textsuperscript{2}{\small{Department of Mathematics, CUNEF University, Madrid, Spain.}}\\%
\textsuperscript{3}{\small{Basque Center for Applied Mathematics (BCAM),  Bilbao, Spain.}}\\%
\textsuperscript{4}{\small{Ikerbasque (Basque Foundation For Sciences), Bilbao, Spain.}}}%
\date{}
\maketitle
\begin{abstract}
Developing efficient methods for solving parametric partial differential equations is crucial for addressing inverse problems.
This work introduces a Least-Squares-based Neural Network (LS-Net) method for solving linear parametric PDEs. It utilizes a separated representation form for the parametric PDE solution via a deep neural network and a least-squares solver. In this approach, the output of the deep neural network consists of a vector-valued function, interpreted as basis functions for the parametric solution space, and the least-squares solver determines the optimal solution within the constructed solution space for each given parameter. The LS-Net method requires a quadratic loss function for the least-squares solver to find optimal solutions given the set of basis functions. In this study, we consider loss functions derived from the Deep Fourier Residual and Physics-Informed Neural Networks approaches. We also provide theoretical results similar to the Universal Approximation Theorem, stating that there exists a sufficiently large neural network that can theoretically approximate solutions of parametric PDEs with the desired accuracy. We illustrate the LS-net method by solving one- and two-dimensional problems. Numerical results clearly demonstrate the method's ability to approximate parametric solutions.
\end{abstract}


\section{Introduction}\label{sec_Introduction}
\noindent Partial Differential Equations (PDEs) play a fundamental role in various fields of science, engineering, and mathematics due to their ability to describe complex phenomena. Parametric PDEs involve problem-dependent parameters -- often related to locally varying, material-dependent properties. Establishing an efficient method for solving parametric PDEs is crucial for tackling inverse problems, where the goal is to determine unknown parameters or functions from observed data. Such inverse problems arise in various fields of science, including medical imaging \cite{bonilla2008inverse}, geophysics \cite{zhdanov2002geophysical}, and electromagnetics  \cite{alvarez2013inversion}. Here, we focus on linear parametric PDEs.


Multiple classical algorithms exist for solving parametric PDEs (e.g., \cite{haasdonk2008reduced, khoromskij2011tensor}). Among them, we highlight the Proper Generalized Decomposition (PGD) technique \cite{ chinesta2011overview, chinesta2011short, garcia2017monitoring, sibileau2018explicit}. The PGD  method is an iterative scheme that approximates the parametric PDE solution using a separated (low-rank) representation of the form
\begin{align}
    u^p(x) := \mathbf{c} (p) \cdot \mathbf{u}(x) = \sum_{n=1}^N c_n(p) u_n(x)
    , \qquad \forall p \in \mathbb{P}
    \label{eq:PGD_decomposition}
\end{align} 
where $p$ is a parameter in the parameter space $\mathbb{P}$, $u^p$ is the approximate parametric solution, and {$\mathbf{c}(p) := (c_1(p), c_2(p), \ldots, c_N(p))$} and {$\mathbf{u}(x) := (u_1(x), u_2(x), \ldots , u_N(x))$} are functions to be determined. The PGD method aims to first determine the basis functions $u_n$ and later compute the coefficients  $c_n$ using various numerical techniques, such as finite elements \cite{discacciati2024overlapping} and Petrov-Galerkin  \cite{nouy2010priori} methods. Despite the numerous advantages of the PGD method for solving parametric PDEs such as computational efficiency and suitability for high-dimensional problems, it still faces significant limitations regarding implementation complexity, including mesh generation and the manual decomposition of the problem into suitable modes.


In the last decade, researchers have explored the use of  Neural Networks (NNs) for solving parametric PDEs \cite{ bhattacharya2021model, brevis2024learning, dal2020data, geist2021numerical, han2018solving, khara2024neufenet,khoo2021solving, kutyniok2022theoretical, uriarte2022finite}. To this end, two different approaches exist: (a) to solve a problem for multiple parameter values using a classical technique (or an NN) and then interpolate it with an NN \cite{bhattacharya2021model, lu2022comprehensive, nelsen2021random, bachmayr2024variationally}, or (b) to directly employ an NN that approximates the full parametric PDE operator (see, e.g., \cite{khoo2019switchnet, kovachki2023neural, li2020fourier, li2020neural, lu2019deeponet, lu2021learning}). Herein, we focus on the approach (b).


 According to the universal approximation theorem \cite{pinkus1999approximation} and subsequent theoretical investigations \cite{chen1995approximation, chen1995universal}, NNs can serve as a potent tool for approximating operators such as solution operators for parametric PDEs. One of the most popular methods in this area is the Neural Operator method \cite{kovachki2023neural}, which learns mappings between infinite-dimensional function spaces. The Neural Operator architecture is designed to be multi-layered, where each layer is composed of linear integral operators and non-linear activation functions. Since the layers are operators stacked together in an end-to-end composition, the overall architecture remains a nonlinear operator. The multi-layered architecture preserves the property of discretization invariance, i.e., the NN structure is independent of the discretization of the underlying function space. The output of the Neural Operator is passed through a pointwise projection operator to transform the output into a function defined on the target domain. This architecture can be instantiated with different practical methods, such as graph-based operators \cite{li2020neural}, low-rank operators \cite{khoo2019switchnet}, and Fourier operators \cite{li2020fourier}. These different instantiations provide flexibility in modeling various types of data and applications. This method also provides the advantages of universal approximation and the ability to handle data at different mesh sizes. The Deep Operator Network (DeepONet) method \cite{lu2019deeponet, lu2021learning}, which maps the values of the input function at a fixed, finite number of sampling points into infinite-dimensional spaces, can be seen as a special case of Neural Operator architectures when restricted to fixed input grids, although one loses the desirable discretization invariance property.

The DeepONet method was originally proposed for learning operators but, in particular, it can be employed for solving parametric PDEs \cite{wang2021learning}.
This method utilizes two NNs: the branch and trunk networks. The branch NN depends upon the parameter $p$, while the trunk NN depends upon the input coordinate $x$. Both NNs are trained in parallel and the approximate solution is obtained as the inner product of the output vectors of the two NNs. Thus, the approximate solution constructed by DeepONet can also be considered as a separated representation of the form \eqref{eq:PGD_decomposition}. 
The main strength of the DeepONet methods is that after constructing the approximate solution \eqref{eq:PGD_decomposition}, solving the parametric PDE for a specific parameter value becomes computationally inexpensive as it only involves forward evaluating two NNs followed by an inner product.
 The Universal Approximation Theorem for operators ensures that the desired accuracy can be achieved under appropriate technical hypotheses via this method using a sufficiently large NN. Additionally, the authors in \cite{lu2021learning} provide some theoretical analysis on the required number of input data to achieve the desired accuracy in learning nonlinear dynamic systems.  This method has gained widespread popularity in multiple applications \cite{goswami2022physics, kovachki2023neural, shukla2024deep, wang2023long}. However, DeepONet typically requires large training datasets of paired input-output observations, which can be expensive to produce. Moreover, its convergence theory does not account for optimization or generalization errors, which may be large. For example, even if the DeepONet is well trained over many instances of PDEs, it may fail to produce an accurate solution if the right-hand side of the PDE is multiplied by a large constant, even though, by linearity, the exact solution will merely be a rescaling of a solution that the network can approximate. By considering the separated representation \eqref{eq:PGD_decomposition}, we see that the failure arises due to poor generalization properties of the functions $c_n$, even if the span of basis functions $u_n$ contains a good approximation.



{In this work, we propose the LS-Net method, which employs the separated representation \eqref{eq:PGD_decomposition}. It utilizes solely one deep NN to construct a vectorial function, whose components are interpreted as the basis functions $u_n(x)$. Then, for each parameter $p$, a least-squares (LS) solver is employed on the solution space spanned by the basis functions to obtain the coefficients $c_n(p)$.  This approach shares similarities with Reduced Order Modeling (ROM) techniques, which are widely utilized in computational science and engineering to enhance computational efficiency. Indeed, both the ROM and the LS-Net methods rely on the fundamental principle of identifying a low-dimensional subspace that effectively approximates the solution manifold. In ROM, this subspace is typically constructed using techniques such as Proper Orthogonal Decomposition \cite{azeez2001proper},  Reduced Basis \cite{rozza2008reduced}, and the PGD \cite{chinesta2011short} methods. In contrast, the LS-Net utilizes an NN to construct these basis functions, thereby identifying a low-dimensional subspace capable of capturing a wide range of solutions to the parametric problem. We provide theoretical results stating that a sufficiently large NN with an LS solver can theoretically approximate solutions within a given desired accuracy. This method also preserves the property of discretization invariance, meaning that whilst appropriate discretization of the loss must be employed to obtain coefficients via the LS system, the obtained basis functions may be utilized with an alternative discretization of the loss. In particular, a finer discretization may be employed when more precise solutions are required. }


On the other side, the LS-Net method also has some limitations, i.e., evaluating the solution for a new parameter $p$ requires the construction and the solution of an LS system, whose computational cost has undesirable scaling properties with respect to the refinement of the discretization. However, this issue can be mitigated once the network is trained, as the LS system can then be constructed parametrically, leading to reduced computational costs. Indeed, this approach eliminates the need for costly integrations at every iteration, making the construction significantly more efficient after performing a single, high-quality integration in advance. In addition, the LS solver requires a quadratic loss function, although this is satisfied in the case of linear PDEs by most NN-based PDE solvers, including Deep Ritz  \cite{yu2018deep}, Double Deep Ritz \cite{uriarte2023deep}, Physics-Informed Neural Networks (PINNs)  \cite{raissi2019physics}, Variational PINNs (VPINNs) \cite{kharazmi2019variational}, Robust VPINNs (RVPINNs) \cite{rojas2024robust}, and Deep Fourier Residual (DFR) \cite{taylor2023deep}.


In this work, we consider two distinct loss functions:  one derived from the DFR method \cite{taylor2023deep}, and another from the PINN method \cite{raissi2019physics}. The PINN loss function is defined based on the strong-form residual of the PDE. {In contrast,  the DFR loss function, a specific case of a VPINN \cite{kharazmi2019variational}, is defined based on the weak-form residual of the PDE.} We use these loss functions to illustrate the strengths and limitations of the proposed parametric PDE solver. 
 

The rest of the paper is organized as follows. 
In Section \ref{Methodology}, we outline the methodology of the LS-Net approach, providing a detailed explanation of parametric residual minimization, the LS-Net solution, loss function discretization, and the neural network structure framework.
Section \ref{Model Problems} provides three numerical examples to verify the theory of the LS-Net method: (a) a damped harmonic oscillator, (b) a one-dimensional Helmholtz equation with impedance boundary conditions, and (c) a two-dimensional transmission problem.
Next, Section \eqref{Conclusion} is dedicated to the conclusions and outlines potential directions for future research. Finally, the Appendix presents the theoretical results concerning the LS-Net method's capabilities in approximating parametric solutions.

\section{Methodology}\label{Methodology}
{This section contains the theoretical foundation and computational framework for solving parametric PDEs using the LS-Net approach. Beginning with the variational formulation of parametric problems, we establish the conditions for well-posedness and introduce the parametric residual minimization framework. This is followed by constructing an NN to approximate the solution space and a discussion on the discretization techniques for implementation. These concepts are outlined in the three subsections.}
\subsection{Parametric residual minimization}
Let $\mathbb{U}$ (trial) and $\mathbb{V}$ (test) denote two Hilbert spaces with corresponding norms $\Vert\cdot\Vert_\mathbb{U}$ and $\Vert\cdot\Vert_\mathbb{V}$, respectively. We consider a parameter space $\mathbb{P}$, and for each $p\in \mathbb{P}$, consider a bounded linear operator $\mathscr{B}^p:\mathbb{U}  \to \mathbb{V}^*$ and $l^p\in \mathbb{V}^*$, where $\mathbb{V}^*$ stands for the topological dual of $\mathbb{V}$. Then, solving a parametric PDE in variational form can be read as: for each $p\in\mathbb{P}$, find $\mathfrak{u}^p\in\mathbb{U}$ such that
\begin{equation}\label{eq_linear weak}
\langle\mathscr{B}^p\mathfrak{u}^p,v\rangle_{\mathbb{V}^* \times \mathbb{V}} = \langle l^p, v\rangle_{\mathbb{V}^* \times \mathbb{V}}, \quad \forall v \in \mathbb{V},
\end{equation}
where $\langle \cdot, \cdot\rangle_{\mathbb{V}^* \times \mathbb{V}}$ is the duality pairing, and $\mathfrak{u}^p$ is the exact solution. As a technical assumption, we take $\mathbb{P}$ to be a compact topological space and assume $\mathfrak{u}^p$ to depend continuously on $p$. {For instance, in this framework, the transmission problem is represented by the operator $\mathscr{B}^p: H_0^1(\Omega) \to H^{-1}(\Omega):=[H_0^1(\Omega)]^*$ given by
\begin{equation*}
\langle\mathscr{B}^p\mathfrak{u}^p,v\rangle_{H^{-1} \times H^1} = \int_\Omega p\nabla \mathfrak{u}^p\cdot\nabla v \, dx,\qquad \forall v\in H^1_0(\Omega),
\end{equation*} 
where $\Omega$ is the domain. } To guarantee the well-posedness of \eqref{eq_linear weak}, we consider the following conditions on $\mathscr{B}^p$ aligned with the hypotheses of the {Babu{\v s}ka}-Lax-Milgram theorem {\cite{babuvska1971error, chen2014inf}}:
\begin{itemize}
\item (Global boundedness) There exists a parameter-independent constant $0<\vartheta < \infty$ such that 
\begin{equation}\label{eq_Continuity condition}
\sup _{u  \in \mathbb{U} \setminus\{0\}} \sup _{v \in \mathbb{V} \setminus\{0\}} \frac{| \langle\mathscr{B}^p u, v\rangle_{\mathbb{V}^* \times \mathbb{V}}|}{\Vert u\Vert_\mathbb{U} \Vert v\Vert_\mathbb{V}} \leq \vartheta, \qquad \forall p \in \mathbb{P}.
\end{equation}
\item (Global weak coercivity) There exists a parameter-independent constant $0<\gamma\leq \vartheta$ such that
\begin{equation}\label{eq_Inf-sup stability condition}
\inf _{u  \in \mathbb{U} \setminus\{0\}} \sup _{v \in \mathbb{V} \setminus\{0\}}\frac{|\langle\mathscr{B}^p u, v\rangle_{\mathbb{V}^* \times \mathbb{V}}|}{ \Vert u\Vert_\mathbb{U} \Vert v\Vert_\mathbb{V}} \geq \gamma, \qquad \forall p \in \mathbb{P}.
\end{equation}
\end{itemize}
As a result, for each $p \in \mathbb{P}$ and any $l^p \in \mathscr{B}^p \mathbb{U}$, problem \eqref{eq_linear weak} admits a unique solution satisfying the following robustness relation between the residual and the error for each $p\in\mathbb{P}$
\begin{equation}\label{eqEquiv}
\frac{1}{\vartheta} \Vert \underbrace{\mathscr{B}^p u - l^p}_{\text{residual}}\Vert_{ \mathbb{V}^*}\leq \Vert \underbrace{u-\mathfrak{u}^p}_{\text{error}}\Vert_{ \mathbb{U}}\leq \frac{1}{\gamma} \Vert \underbrace{\mathscr{B}^p u - l^p}_{\text{residual}} \Vert_{ \mathbb{V}^*},\qquad u\in\mathbb{U},
\end{equation} 
where $\displaystyle \Vert\cdot\Vert_{\mathbb{V}^*} = \sup_{v\in\mathbb{V}\setminus\{0\}} \frac{|\langle\cdot,v\rangle_{\mathbb{V}^* \times \mathbb{V}}|}{\Vert v\Vert_{\mathbb{V}}}$.
This implies that the unique solution of \eqref{eq_linear weak} is given by
\begin{equation}\label{eq_exact solution}
 \mathfrak{u}^p =\arg\min_{u\in \mathbb{U}} \Vert {\mathscr{B}^p u - l^p}\Vert_ {\mathbb{V}^*}^2,\qquad p\in\mathbb{P}.
\end{equation}
\begin{remark}
Spaces $\mathbb{U}$ and $\mathbb{V}$ can be generalized to  parameter-dependent spaces $\mathbb{U}^p$ and $\mathbb{V}^p$ for $p \in \mathbb{P}$, respectively. 
Parameter-dependent spaces allow the framework to include more appropriate norms that may provide more desirable constants $\vartheta$ and $\gamma$, as in \cite{bachmayr2024variationally}. {For instance, for the transmission problem with Dirichet BCs, the parameter-depended energy norm is defined as
\begin{equation}
\Vert u^p \Vert_{\mathbb{U}^p} = \int p(\nabla u^p)^2\, dx, \qquad \forall u^p \in \mathbb{U}^p.
\end{equation}}
This flexibility is crucial for addressing indefinite, non-symmetric, or singularly perturbed problems. Additionally, for certain complex PDEs, the stability and well-posedness of the variational formulation \eqref{eq_linear weak} might only be achievable by allowing the trial and test spaces to depend on the parameter $p$. Herein, we simplify the problem and theoretical analysis by considering parameter-independent spaces.
\end{remark}


\subsection{LS-Net method}
{Let $\alpha$ denote the set of learnable parameters of an NN yielding the output vector}
\begin{equation}
    x\mapsto\mathbf{u}^\alpha(x) := (u^\alpha _1(x), u^\alpha _2(x), \ldots , u^\alpha _N(x)).
\end{equation} 
{By considering the approximate solution space $\mathbb{U}^\alpha = \text{span}\{u^\alpha_n\}_{n=1}^N$, we represent the approximate parametric solution of \eqref{eq_linear weak} in the form of }
\begin{align}
     \mathfrak{u}^p(x) \approx u^{p, \alpha}(x) := \mathbf{c}^{p,\alpha} \cdot \mathbf{u}^\alpha(x),\qquad p\in\mathbb{P},
    \label{eq:PGD_decomposition_NN}
\end{align}
where $\mathbf{c}^{p,\alpha}:=(c^{p,\alpha}_{1},c^{p,\alpha}_{2}, \ldots, c^{p,\alpha}_{N})$ is the vector of optimal coefficients meeting the LS condition in \eqref{eq_exact solution}, i.e.
\begin{equation}\label{eq_optimal_coeff}
    \mathbf{c}^{p,\alpha} = \arg\min_{\mathbf{c}} \Vert {\mathscr{B}^p(\mathbf{c}\cdot\mathbf{u}^\alpha) - l^p}\Vert_ {\mathbb{V}^*}^2.
\end{equation} 
As a result, $u^{p, \alpha}$ is the minimum-residual element lying on {$\mathbb{U}^\alpha$}. 
For non-parametric problems (i.e., for fixed parameter $p \in \mathbb{P}$), \cite{uriarte2024optimizing} proposed an optimization scheme that trains $\alpha$ (adapts $\mathbb{U}^\alpha$) based on the optimal coefficients to minimize the following loss function \eqref{eq_optimal_coeff}
\begin{equation}\label{eq_training_nonparametric}
    \mathcal{L}^p(\alpha)\;:= \; \min_{\mathbf{c}} \Vert {\mathscr{B}^p(\mathbf{c}\cdot\mathbf{u}^\alpha) - l^p}\Vert_ {\mathbb{V}^*}^2 =  \Vert {\mathscr{B}^p u^{p,\alpha} - l^p}\Vert_ {\mathbb{V}^*}^2.
\end{equation} 

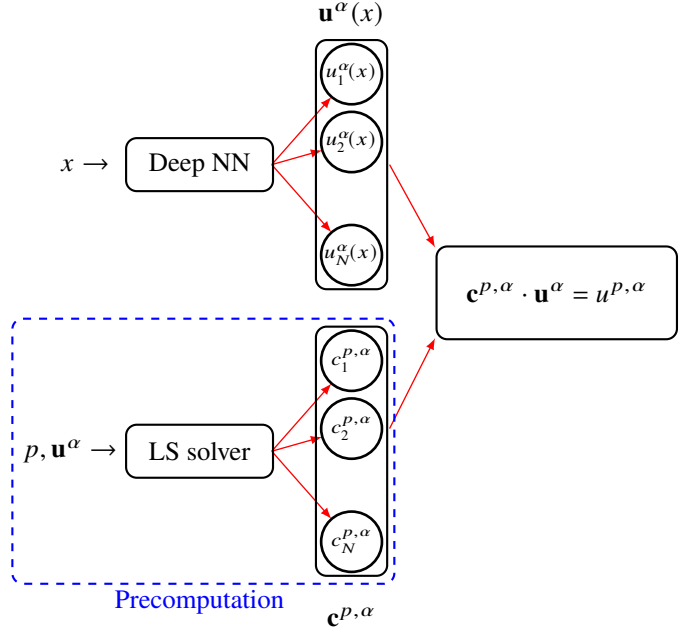
\begin{figure}[htbp]
\begin{center}
\scalebox{1}{\usetikzlibrary{decorations.pathmorphing} 
\begin{tikzpicture}[
punkt/.style={
   rectangle,
   rounded corners,
   draw=black, thick,
   text width=2em,
   minimum height=1em,
   text centered}
]

\draw (-2,1.5) node{ $x \rightarrow$};
\node[punkt, minimum width = 5.5em, minimum height = 2em] (block) at (-0.5, 1.5) {};

\draw (-0.5, 1.5) node{Deep NN};


\Vertex[x=1.5,y=2.7,label=$u_1^\alpha\hspace{-0.05cm}(x)$,color=red!20!white,size=0.8,fontsize=\scriptsize]{A11}
\Vertex[x=1.5,y=1.8, label=$u_2^\alpha\hspace{-0.05cm}(x)$,color=red!20!white,size=0.8,fontsize=\scriptsize]{A12}
\Vertex[x=1.5,y=0.3, label=$u_N^\alpha\hspace{-0.05cm}(x)$,color=red!20!white,size=0.8,fontsize=\scriptsize]{A13}
\draw (1.5,1.2) node{ $\vdots$};

\Edge[color=black, Direct,lw=0.5pt](0.5, 1.5)(A11)
\Edge[color=black, Direct,lw=0.5pt](0.5, 1.5)(A12)
\Edge[color=black, Direct,lw=0.5pt](0.5, 1.5)(A13)

\draw (1.5,3.5) node{ ${\color{red}\mathbf{u}^\alpha} (x)$};

\node[punkt, minimum width = 14.3em, minimum height = 12em, draw, black] (block) at (-0.15, 1.8) {};

\draw (4.6,1.5) node{ $p, {\color{red}\mathbf{u}^\alpha} \rightarrow$};
\node[punkt, minimum width = 12.5em, minimum height = 4em] (block) at (7.6, 1.5) {};

\draw[thick, ->] (2.5, 1.7) -- (3.5, 1.7); %

\draw (7.6, 1.85) node{LS solver};

\draw (7.6, 1.15) node{$
    \displaystyle \min_{\mathbf{c}} \Vert \mathscr{B}^p(\mathbf{c}\cdot {\color{red}\mathbf{u}^\alpha}) - l^p\Vert_ {\mathbb{V}^*}^2$};


\Vertex[x=10.8,y=2.7,label=$c^{p, \alpha}_{1}$,color=blue!20!white,size=0.8,fontsize=\scriptsize]{A11}
\Vertex[x=10.8,y=1.8, label=$c^{p, \alpha}_{2}$,color=blue!20!white,size=0.8,fontsize=\scriptsize]{A12}
\Vertex[x=10.8,y=0.3, label=$c^{p, \alpha}_{N}$,color=blue!20!white,size=0.8,fontsize=\scriptsize]{A13}
\draw (10.8,1.2) node{ $\vdots$};

\Edge[color=black, Direct,lw=0.5pt](9.8, 1.5)(A11)
\Edge[color=black, Direct,lw=0.5pt](9.8, 1.5)(A12)
\Edge[color=black, Direct,lw=0.5pt](9.8, 1.5)(A13)

\draw (10.8,3.5) node{ ${\color{blue}\mathbf{c}^{p,\alpha}}$};
\node[punkt, minimum width = 22.7em, minimum height = 12em, draw, black] (block) at (7.6, 1.8) {};


\draw[thick] (-0.2, -0.3) -- (-0.2, -1.2); 
\draw[thick, ->] (-0.2, -1.2) -- (2.2, -1.2);  

\draw[thick] (7.6, -0.3) -- (7.6, -1.2); 
\draw[thick, ->] (7.6, -1.2) -- (6.2, -1.2); 

\node[punkt, minimum width = 10em, minimum height = 2.5em] (block) at (4.2, -1.2) {};
\draw (4.2, -1.2) node{$\displaystyle u^{p, \alpha}(x) = {\color{blue}\mathbf{c}^{p,\alpha}} \cdot {\color{red}\mathbf{u}^\alpha}(x)$};

\end{tikzpicture}}
\end{center}		
\caption{Schematic of the LS-Net architectures for solving parametric PDEs.}\label{fig_ARC}
\end{figure}

In this work, we propose the natural extension of the above idea to our parametric problem \eqref{eq_linear weak}, training $\alpha$ to ensure that the solution space $\mathbb{U}^\alpha$  is suited for all $p\in\mathbb{P}$. To this end, we consider the parametric loss function given by
\begin{equation}\label{eq_continuum_Loss}
\mathcal{L}_\mu(\alpha) := \int_\mathbb{P}\mathcal{L}^p(\alpha) \,d\mu(p),
\end{equation} 
where $\mu$ denotes a probability measure on $\mathbb{P}$. This allows us to train $\alpha$ so that the LS solutions \eqref{eq_optimal_coeff} yield accurate results across the entire parameter space. Note that when considering $\mathbb{P}=\{p\}$ in \eqref{eq_continuum_Loss}, we recover the non-parametric formulation \eqref{eq_training_nonparametric}. Figure \ref{fig_ARC} presents a schematic of the LS-Net architecture.

Below, we state the universal approximation theorem for this parametric framework. Its proof is deferred to Appendix A.

\begin{theorem}
Consider conditions \eqref{eq_linear weak}-\eqref{eq_Inf-sup stability condition} with $\mathbb{U} = \mathbb{V} = H^k(\Omega)$ where  {$H^k$ is the standard Sobolev space of order} $k\in\mathbb{N}$ {and for $N_0 \in\mathbb{N}$, $\Omega\subset \mathbb{R}^{N_0} $ is a bounded and Lipschitz domain.}  Moreover, assume that $\mathbb{P}\ni p\mapsto \mathfrak{u}^p\in\mathbb{U}$ is $\mu$-measurable and continuous. Then, for every $\epsilon>0$, there exists $\alpha_\epsilon$ and $N\in\mathbb{N}$ yielding a fully-connected feed-forward neural network $\mathbf{u}^{\alpha_\epsilon}:\Omega\longrightarrow\mathbb{R}^N$ such that
\begin{equation}
\int_\mathbb{P}||u^{p, \alpha_\epsilon}-\mathfrak{u}^p||_{\mathbb{U}}^2\,d\mu(p)<\epsilon,\qquad
\int_\mathbb{P} \Vert \mathscr{B}^p u^{p, \alpha_\epsilon} - l^p \Vert _{\mathbb{V}^*}^2\,d\mu(p) <\epsilon,
\end{equation}
where $u^{p,\alpha_\epsilon} = \mathbf{c}^{p,\alpha_\epsilon}\cdot\mathbf{u}^{\alpha_\epsilon}:\Omega\longrightarrow\mathbb{R}$ with $\mathbf{c}^{p,\alpha_\epsilon}$ as in \eqref{eq_optimal_coeff}.
\end{theorem}

\subsection{Loss function discretization}
To enable a practical implementation, the continuum loss function \eqref{eq_continuum_Loss} must be discretized, which involves discretizing the integral over the parameter space $\mathbb{P}$, followed by addressing the challenges associated with evaluating the dual norm and the optimal coefficients $\mathbf{c}^{p,\alpha}$.


\subsubsection*{Integral approximation over the parameter space}\label{Parameter space discretization}
\noindent To discretize \eqref{eq_continuum_Loss} over the parameter space $\mathbb{P}$, we employ Monte Carlo integration, resulting in
\begin{equation}\label{eq_parameter-discretization_loss}
\mathcal{L}_\mu (\alpha) \approx \frac{1}{|P|} \sum \limits_{p\in P} \mathcal{L}^p(\alpha), 
\end{equation}
where $P\subset\mathbb{P}$ denotes a batch of independent and identically distributed samples following the probability measure $\mu$.

\subsubsection*{The dual norm discretization}\label{Approximation of the dual norm}
To compute the semi-discretized loss function \eqref{eq_parameter-discretization_loss}, it is essential to approximate the dual norm $\Vert \cdot \Vert_{\mathbb{V}^*}$ twice: first, to determine the optimal coefficient vector $\mathbf{c}^{p, \alpha}$, and subsequently, to compute the loss value. In this context, we will examine two methods:  the PINN and the DFR method. 

\begin{subequations}\label{PINNandDFR}
\begin{itemize}
    \item \textbf{PINN method:} Consider $\mathbb{V} = L^2$. Since functionals in $L^2$ are easily identifiable with elements of $L^2$ itself, the discretization of the dual norm can be reduced to applying a numerical quadrature rule. For simplicity, we employ the midpoint quadrature rule. Indeed, for $J \in \mathbb{N}$ and quadrature points $\{ x_j\}_{j=1}^{J}$, we have
\begin{equation}\label{eq_PINN method}
    \Vert \mathscr{B}^{p}(\mathbf{c}\cdot \mathbf{u}^\alpha) - l^{p} \Vert _{L^2}^2 \approx \frac{{\text{Vol}(\Omega)}}{J} \sum_{j=1}^{J} (\mathscr{B}^{p}(\mathbf{c}\cdot \mathbf{u}^\alpha)(x_j) - l^{p}(x_j))^2 = \Vert \mathbf{B}^{p,\alpha}\mathbf{c}-\mathbf{l}^{p}\Vert^2_2, \qquad p \in \mathbb{P},
    \end{equation}
    where $(\mathbf{B}^{p,\alpha})_{jn}={\sqrt{\frac{\text{Vol}(\Omega)}{J}}}(\mathscr{B}^p u^\alpha_n)(x_j)$ is a $J\times N$ matrix, $(\mathbf{l}^{p})_{j}={\sqrt{\frac{\text{Vol}(\Omega)}{J}}}l^p(x_j)$ is a vector of size $J$, and $\Vert\cdot\Vert_2$ is the Euclidean norm. PINNs are one of the most appealing methods due to their simplicity in implementing autodiff algorithms for optimization and efficient training \cite{raissi2017physics,raissi2019physics}. 
    
    \item \textbf{DFR method:} 
    Following Parseval's identity in separable Hilbert spaces, we can express $\Vert \cdot \Vert _{\mathbb{V}^*}$ as a series expansion consisting of duality pairings with an orthonormal basis $\{v_m\}_{m=1}^\infty$ of the test space $\mathbb{V}$, i.e.
    
    \begin{equation}\label{eq_Dualnorm}
    \Vert \mathscr{B}^{p}(\mathbf{c}\cdot \mathbf{u}^\alpha) - l^{p} \Vert _{\mathbb{V}^*}^2= \sum\limits_{m=1}^\infty \langle\mathscr{B}^{p}(\mathbf{c}\cdot \mathbf{u}^\alpha) - l^p, v_m\rangle_{\mathbb{V}^* \times \mathbb{V}}^2, \qquad p \in \mathbb{P},
    \end{equation}
    where $p \in P$. A natural discretization of the dual norm is achieved by truncating the series expansion \eqref{eq_Dualnorm}. Therefore, by selecting $M \in \mathbb{N}$, we have
    \begin{equation}\label{eq_Dualnorm_0}
     \Vert \mathscr{B}^{p}(\mathbf{c}\cdot \mathbf{u}^\alpha) - l^{p} \Vert _{\mathbb{V}^*}^2 \approx \sum\limits_{m=1}^M \langle\mathscr{B}^{p}(\mathbf{c}\cdot \mathbf{u}^\alpha) - l^p, v_m\rangle_{\mathbb{V}^* \times \mathbb{V}}^2 = \Vert \mathbf{B}^{p,\alpha} \mathbf{c} - \mathbf{l}^{p} \Vert_2^2,\qquad p\in P,
    \end{equation}
    where $(\mathbf{B}^{p,\alpha})_{mn}=\langle\mathscr{B}^p u^\alpha_n, v_m\rangle_{\mathbb{V}^*\times\mathbb{V}}$ is a $M\times N$ matrix and $(\mathbf{l}^{p})_{j}=\langle l^p, v_m\rangle_{\mathbb{V}^*\times\mathbb{V}}$ is a vector of size $M$.
    The DFR method employs the orthogonal eigenvectors (in $H^1$) of the weak-form Laplacian as test basis functions $\{v_m\}_{m=1}^M$. The eigenvectors of the Laplacian are naturally ordered by taking the eigenvalues to be increasing, with the truncation acting like a low-pass filter on the residual, capturing its dominant, low-frequency behavior. In Cartesian products of intervals, the eigenfunctions can be expressed via sines and cosines of varying frequencies, and we refer to \cite[Table 1]{taylor2023deep} for further details on these basis functions. Note that the duality pairings in \eqref{eq_Dualnorm_0} are, in practice, approximated using an appropriate quadrature rule, for which we again select the midpoint rule for simplicity.

\end{itemize}
\end{subequations}

\noindent 
{The computation of the loss function is outlined in Algorithm \ref{Algorithm_1}. }

\begin{algorithm}[H]\label{Algorithm_1}
{
\begingroup

\SetAlgoLined
Sample a finite batch of parameter values $P$.\\
Construct $\{\mathbf{B}^{p,\alpha},\mathbf{l}^p\}_{p\in P}$.\\
Compute all vectors of optimal coefficients $\{\mathbf{c}^{p,\alpha}\}_{p\in P}$, by solving the following LS problems
{\begin{equation}\label{eq_optimal_coeff_al}
\displaystyle \mathbf{c}^{p,\alpha} =  \min_{\mathbf{c}} \Vert\mathbf{B}^{p,\alpha}\mathbf{c}-\mathbf{l}^p\Vert_2^2, \qquad p \in P.
\end{equation}}\\
Estimate the resulting minimum-residual average 
\begin{equation}\label{eq_Loss-D}
\mathcal{L}_\mu(\alpha)\approx \frac{1}{\vert P\vert}\sum_{p\in P} \Vert\mathbf{B}^{p,\alpha}\mathbf{c}^{p,\alpha}-\mathbf{l}^p\Vert_2^2.
\end{equation} 
\caption{Implementation for computing the loss function}\label{train_step_guidelines}
\endgroup
}
\end{algorithm}

\noindent {
In Algorithm \ref{Algorithm_1}, $\mathbf{c}^{p,\alpha}$ is indeed an approximation to the continuum-level coefficients indicated in \eqref{eq_optimal_coeff}. However, for simplicity, we retain the same notation to refer to this approximate version. 
To solve the LS problem \eqref{eq_optimal_coeff_al}, we utilized TensorFlow's LS solver, which typically resolves the normal equations using Cholesky decomposition.} 
The matrix $\mathbf{B}^{p,\alpha}$ can be formulated as a combination of specific parameter-independent matrices, allowing to reformulate the LS system parametrically. Therefore, instead of performing costly integrations at each iterative step, only a single high-precision integration is required, leading to a significant reduction in overall computational cost. 

We also emphasize the convenience of using forward-mode automatic differentiation in space instead of its usual reverse mode when evaluating the derivatives of the spanning functions $u^\alpha_n$, as the differences in computational cost are significant \cite[Section 3]{uriarte2024optimizing}. 
For further details on our implementation, see \cite{software_LS-Net4ParametricPDEs}.


\subsection{Neural network architecture}\label{Neural Network structure} 
Let us consider the following fully-connected NN with depth $K$ and the set of learnable parameters $\alpha$

\begin{equation}
\bar{\mathbf{u}}^\alpha = L_{K} \circ L_{K - 1} \circ \ldots  \circ  L_1 \circ L_0, 
\end{equation}
where $\bar{\mathbf{u}}^\alpha = (\bar{u}^\alpha_1, \bar{u}^\alpha_2, \ldots, \bar{u}^\alpha_N)^T$ and $L_0 = x$ are the output and input vectors of the network, respectively. $L_k$ is the layer function and it is defined as follows
\begin{equation}\label{eq_layers}
L_k =\sigma( A_k L_{k-1} + b_k), k = 1, \ldots, K, 
\end{equation}
where $\sigma(\cdot)$ is the activation function that acts component-wise on vectors. Here, we utilize the sigmoid activation function. The matrices {$A_k\in\mathbb R^{N_k} \times \mathbb R^{N_{k-1}}$} and the vectors $b_k\in\mathbb{R}^{N_{k}}$ are called weights and biases, respectively, where their components collectively form the set $\alpha$. In this representation, {$N_0$ is the dimension of the problem}, $N_{K} = N$ is the dimension of the output, and $N_{k}$ for $ k =1, \ldots, K-1$ denotes the number of neurons in the $k$-th layer. {Note that the success of an NN heavily depends on the architecture and proper selection of hyperparameters (e.g., initialization, learning rate), as different initializations can yield distinct results. However, a detailed analysis of these factors is beyond the scope of this work.}

In the context of utilizing NNs for solving PDEs, boundary conditions are enforced through two distinct approaches: (a) they are incorporated weakly and the solution space does not satisfy the boundary conditions, or (b) they are applied strongly to the outputs of the NN, and thus, the solution space enforces the boundary conditions. Here, we explain the implementation of approach (b). For homogeneous {Dirichlet} boundary conditions, we apply a cut-off function $\varphi$ to $\bar{\mathbf{u}}^\alpha$, where $\varphi(x) = 0$ on the segments of the boundary for which the solution space includes the boundary conditions \cite{berrone2023enforcing}. Given that any nonhomogeneous PDE problem can be transformed into a homogeneous one, this methodology can be readily adapted to handle inhomogeneous boundary conditions by applying an appropriate shift function. {Moreover, this approach naturally extends to the case where both Dirichlet and Neumann conditions are imposed at a specific point in one dimension. This is achieved by utilizing a cutoff function with $\phi'(x)=0$ and applying an appropriate lifting function.}

Moreover, for model problems exhibiting discontinuities in the gradient (e.g., \ref{One-dimensional Helmholtz equation} and \ref{Two-dimensional transmission problem}), we employ the ideas of Regularity Conforming NNs \cite{taylor2024regularity}. We introduce a problem dependent function $\psi = (\psi_1, \ldots, \psi_N)^T$ of lower regularity, which multiplies $\bar{\mathbf{u}}^{\alpha}$ component-wise. This function ensures the desired regularity of the solution and facilitates more accurate approximations in the discontinuities.

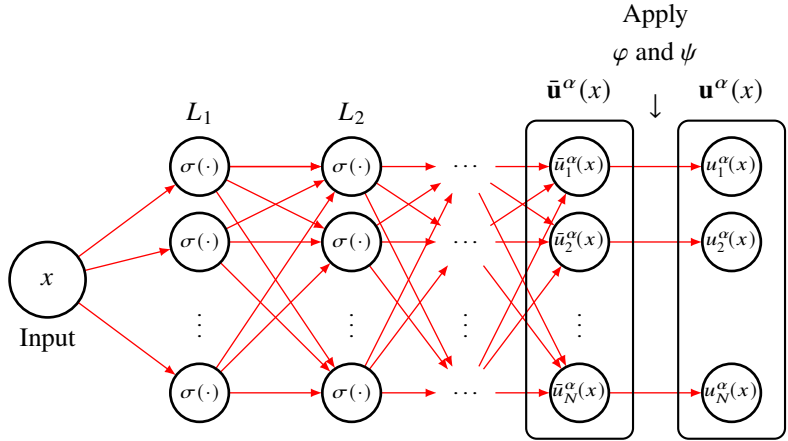
\begin{figure}[htbp]
	\begin{center}
	\begin{tikzpicture}[
punkt/.style={
   rectangle,
   rounded corners,
   draw=black, thick,
   text width=2em,
   minimum height=1em,
   text centered}
]


\Vertex[x=10,y=1.5,style={opacity=0.0,color=white!0},size=0.00000001]{inv1}
\Vertex[x=10,y=0,style={opacity=0.0,color=white!100}]{inv2}

\Vertex[x=0,y=1.5,label=$ x$,color=white,size=1,fontsize=\normalsize]{X}
\draw (0,0.7) node {Input};

\Vertex[x=2,y=3,label=$\sigma (\cdot)$,color=white,size=0.77,fontsize=\scriptsize]{A11}
\Vertex[x=4,y=3,label=$\sigma (\cdot)$,color=white,size=0.77,fontsize=\scriptsize]{A21}
\Vertex[x=2,y=2,label=$\sigma (\cdot)$,color=white,size=0.77,fontsize=\scriptsize]{A12}
\Vertex[x=4,y=2,label=$\sigma (\cdot)$,color=white,size=0.77,fontsize=\scriptsize]{A22}
\Vertex[x=2,y=0,label=$\sigma (\cdot)$,color=white,size=0.77,fontsize=\scriptsize]{A13}
\Vertex[x=4,y=0,label=$\sigma (\cdot)$,color=white,size=0.77,fontsize=\scriptsize]{A23}

\Vertex[x=7,y=3,label=$\bar{u}_{1}^\alpha \hspace{-0.05cm}(x)$ ,color=white,size=0.77,fontsize=\scriptsize]{A31}
\Vertex[x=7,y=2,label=$\bar{u}_{2}^\alpha \hspace{-0.05cm}(x)$ ,color=white,size=0.77,fontsize=\scriptsize]{A32}
\Vertex[x=7,y=0,label=$\bar{u}_{N}^\alpha \hspace{-0.05cm}(x)$ ,color=white,size=0.77,fontsize=\scriptsize]{A33}

\Vertex[x=9,y=3,label=$u_{1}^\alpha \hspace{-0.05cm}(x)$ ,color=white,size=0.77,fontsize=\scriptsize]{A41}
\Vertex[x=9,y=2,label=$u_{2}^\alpha \hspace{-0.05cm}(x)$ ,color=white,size=0.77,fontsize=\scriptsize]{A42}
\Vertex[x=9,y=0,label=$u_{N}^\alpha \hspace{-0.05cm}(x)$ ,color=white,size=0.77,fontsize=\scriptsize]{A43}

\draw (2,3.7) node {$L_1$};
\draw (4,3.7) node {$L_2$};
\draw (7,4) node {$\bar{\mathbf{u}}^\alpha (x)$};
\draw (9,4) node {$\mathbf{u}^\alpha (x)$};

\Vertex[x=5.5,y=3,style={color=white},label=$\color{black}\hdots$,size=0.75]{Ah1}
\Vertex[x=5.5,y=2,style={color=white},label=$\color{black}\hdots$,size=0.75]{Ah2}
\Vertex[x=5.5,y=0,style={color=white},label=$\color{black}\hdots$,size=0.75]{Ah3}

\Edge[color=red, Direct,lw=0.5pt](X)(A11)
\Edge[color=red, Direct,lw=0.5pt](X)(A12)
\Edge[color=red, Direct,lw=0.5pt](X)(A13)

\Edge[color=white,label={$\color{black}\vdots$}](A12)(A13)
\Edge[color=white,label={$\color{black}\vdots$}](A22)(A23)
\Edge[color=white,label={$\color{black}\vdots$}](Ah2)(Ah3)
\Edge[color=white,label={$\color{black}\vdots$}](A32)(A33)

\Edge[color=red, Direct,lw=0.5pt](A11)(A21)
\Edge[color=red, Direct,lw=0.5pt](A11)(A22)
\Edge[color=red, Direct,lw=0.5pt](A11)(A23)
\Edge[color=red, Direct,lw=0.5pt](A11)(A21)
\Edge[color=red, Direct,lw=0.5pt](A12)(A21)
\Edge[color=red, Direct,lw=0.5pt](A12)(A22)
\Edge[color=red, Direct,lw=0.5pt](A12)(A23)
\Edge[color=red, Direct,lw=0.5pt](A13)(A21)
\Edge[color=red, Direct,lw=0.5pt](A13)(A22)
\Edge[color=red, Direct,lw=0.5pt](A13)(A23)

\Edge[color=red, Direct,lw=0.5pt](A23)(Ah1)
\Edge[color=red, Direct,lw=0.5pt](A23)(Ah2)
\Edge[color=red, Direct,lw=0.5pt](A23)(Ah3)
\Edge[color=red, Direct,lw=0.5pt](A21)(Ah1)
\Edge[color=red, Direct,lw=0.5pt](A21)(Ah2)
\Edge[color=red, Direct,lw=0.5pt](A21)(Ah3)
\Edge[color=red, Direct,lw=0.5pt](A22)(Ah1)
\Edge[color=red, Direct,lw=0.5pt](A22)(Ah2)
\Edge[color=red, Direct,lw=0.5pt](A22)(Ah3)

\Edge[color=red, Direct,lw=0.5pt](Ah1)(A31)
\Edge[color=red, Direct,lw=0.5pt](Ah1)(A32)
\Edge[color=red, Direct,lw=0.5pt](Ah1)(A33)
\Edge[color=red, Direct,lw=0.5pt](Ah2)(A31)
\Edge[color=red, Direct,lw=0.5pt](Ah2)(A32)
\Edge[color=red, Direct,lw=0.5pt](Ah2)(A33)
\Edge[color=red, Direct,lw=0.5pt](Ah3)(A31)
\Edge[color=red, Direct,lw=0.5pt](Ah3)(A32)
\Edge[color=red, Direct,lw=0.5pt](Ah3)(A33)

\Edge[color=red, Direct,lw=0.5pt](A31)(A41)
\Edge[color=red, Direct,lw=0.5pt](A32)(A42)
\Edge[color=red, Direct,lw=0.5pt](A33)(A43)

\node[punkt, minimum width = 4em, minimum height = 12em] (block) at (9,1.5) {};
\node[punkt, minimum width = 4em, minimum height = 12em] (block) at (7,1.5) {};

\draw (8,3.8) node{ $ \downarrow $};
\draw (8,4.5) node{$\varphi$ and $\psi$};
\draw (8,5) node{Apply};


%

\end{tikzpicture}
	\end{center}
	\caption{NN architecture.}
	\label{fig_NN_architecture}
\end{figure}

Consequently, the final NN architecture is defined as follows
\begin{equation}\label{eq_regularity-conforming}
\mathbf{u}^\alpha(x) = \varphi (x)  \bar{\mathbf{u}}^\alpha(x) \odot \psi(x),
\end{equation}
where the symbol $\odot$ represents the component-wise multiplication. When we consider problems with smooth solutions, we take $\psi=1$. Figure \ref{fig_NN_architecture} illustrates the overall NN architecture. The structure of $\mathbf{u}^\alpha$ introduces a negligible computational cost compared to $\bar{\mathbf{u}}^\alpha$ while allowing greater approximation capability of the solution. The precise choices of functions $\varphi$ and $\psi$ are problem-dependent and will be described in Section \ref{Model Problems} for their corresponding examples.

\section{Model Problems}\label{Model Problems}

\subsection{Damped harmonic oscillator}
The motion of a mass on a damped spring is a classical Ordinary Differential Equation (ODE) that can show a variety of distinct behaviors, typically described by
\begin{equation}\label{eq_damped problem}
m \mathfrak{u}_{tt}^p + c \mathfrak{u}_t^p  + k \mathfrak{u}^p = 0, \qquad t \in (0, T),
\end{equation}
where $p = (m, c, k)$ is the triplet of parameters and $m$, $c$, and $k$ are positive constants representing the mass, damping constant, and spring constant, respectively. The solution $\mathfrak{u}^p$ is a twice differentiable function of $t$ that represents the one-dimensional displacement of the mass from its equilibrium position in the absence of external forces. For fixed initial conditions $\mathfrak{u}^p(0)=\mathfrak{u}_0^p$ and $\mathfrak{u}^p_t(0)=\mathfrak{v}_0^p$, equation \eqref{eq_damped problem} has a unique, smooth solution. The general solution of Equation \eqref{eq_damped problem} takes the form $A \exp(\lambda_1 t) + B \exp(\lambda_2 t)$, where $\lambda_1$ and $\lambda_2$ are the roots of the characteristic polynomial $m \lambda^2 +c \lambda +k$. We remark that the span of all such solutions when varying just one parameter is infinite dimension.

Equation \eqref{eq_damped problem} is particularly interesting because, despite its simplicity, it produces a wide range of behaviors, depending on the sign of the discriminant $c^2-4mk$: when negative, it produces underdamping, which results in oscillatory motion; when positive, it produces overdamping, characterized by an exponential decay toward zero; when it is equal to zero, its results in critical damping, where the solution rapidly approaches zero without oscillating.

Herein, we use the LS-Net method to approximate the high-dimensional solution space of problem \eqref{eq_damped problem} using a low-dimensional, discretized trial space. Without loss of generality, we simplify equation \eqref{eq_damped problem} by reducing the number of parameters through division by the damping constant $c$. This yields
\begin{equation}\label{eq_damped_problem_2p}
p_1 \mathfrak{u}_{tt}^p + \mathfrak{u}_t^p + p_2 \mathfrak{u}^p = 0, \qquad t \in (0, T],
\end{equation}
where $p_1 = m/c$ and $p_2 = k/c$. Also, we set  $T = 10$ and consider initial conditions {$\mathfrak{u}_0^p = 0$} and $\mathfrak{v}_0^p = -50$. 

{We employ an NN of the form described in \eqref{eq_regularity-conforming} with $K=3$, $N_1 = N_2 = 5$, and $N_3 =40$. This network is regarded as low-dimensional, as its output dimension of 40 is relatively small compared to the infinite-dimensional solution space.}
To enforce the boundary conditions, we employ the cut-off function $\varphi(t) = t^2$. As the solutions are smooth, we do not require a regularity conforming NN. Note that, based on the inhomogeneous initial conditions, the solution is obtained by employing a lift, i.e. $u^{p, \alpha} -50 t $.  Additionally, we employ a discretized loss function based on the PINN approach with the midpoint quadrature rule and $1000$ quadrature points on the interval $(0, T]$. The training and validation datasets each consist of 500 randomly sampled parameter values,  generated for $p_1$ and $p_2$ from a log-uniform and uniform distribution over the range $(10^{-1.5}, 10^{1.5})$, respectively. Note that we set the validation dataset fixed throughout, while the training dataset is resampled in the parameter space randomly at each training epoch. {Moreover, for the validation set in time, we used 1005 fixed quadrature points, distinct from those used during training. 
The considered NN is trained using the Adam optimizer over  $10^{4}$ iterations, with a learning rate that decreases exponentially during the training process. Let $\lambda_0$ and $\lambda_E$ represent the initial and final learning rates, respectively, where $E$ is the total number of iterations ($E = 10^{4}$ in this case). The learning rate at iteration $e$ for $e = 1, \ldots, E-1$  is updated at each iteration using the following formula
\begin{equation}
\lambda_e = \lambda_0 ( \frac{\lambda _E}{\lambda _0})^{\frac{e}{E}}. 
\end{equation}
This formulation ensures a smooth and controlled exponential decay of the learning rate over the total number of epochs. In this example, the initial learning rate is set to $10^{-3}$, and the final learning rate is set to $10^{-4}$.}

\begin{figure}[H]

	\begin{center}
	
	\begin{subfigure}[a]{1\textwidth}
	\begin{center}
	\includegraphics[width=1\textwidth]{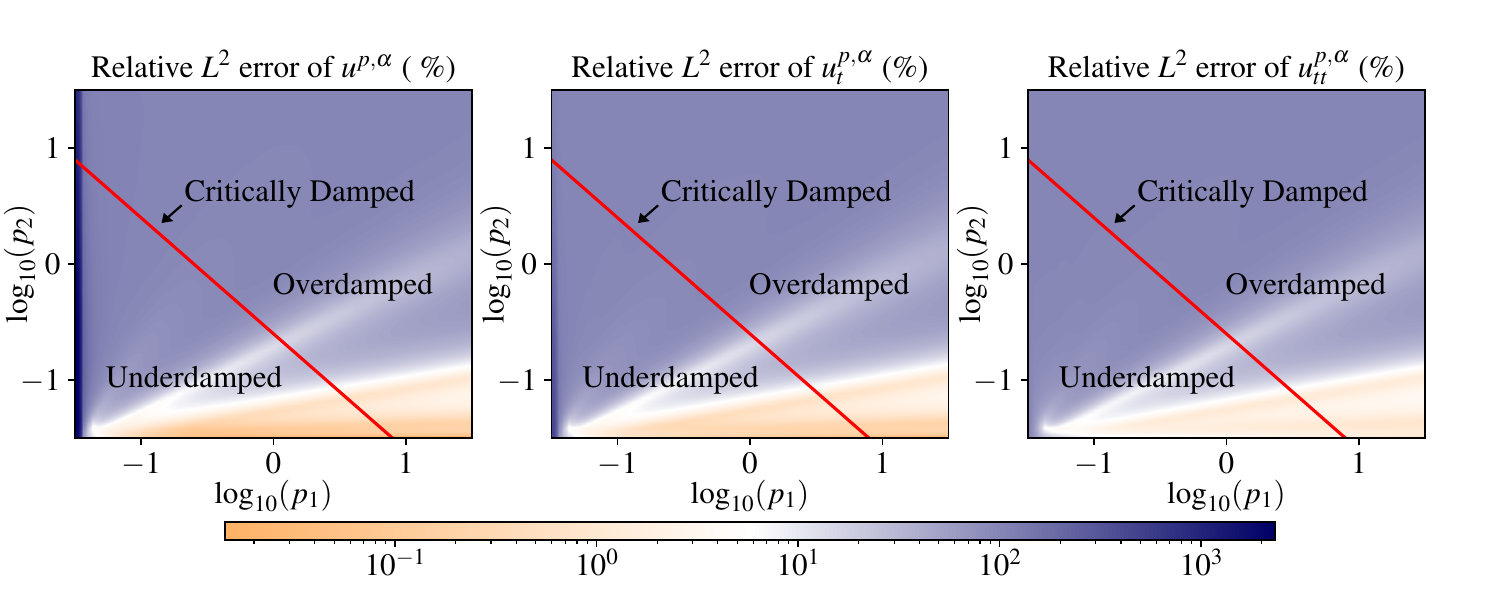}
	\caption{Before training. }\label{fig_PINN_L2_Error_a}
	\end{center}
    \end{subfigure}

	\begin{subfigure}[b]{1\textwidth}
		\begin{center}
	\includegraphics[width=1\textwidth]{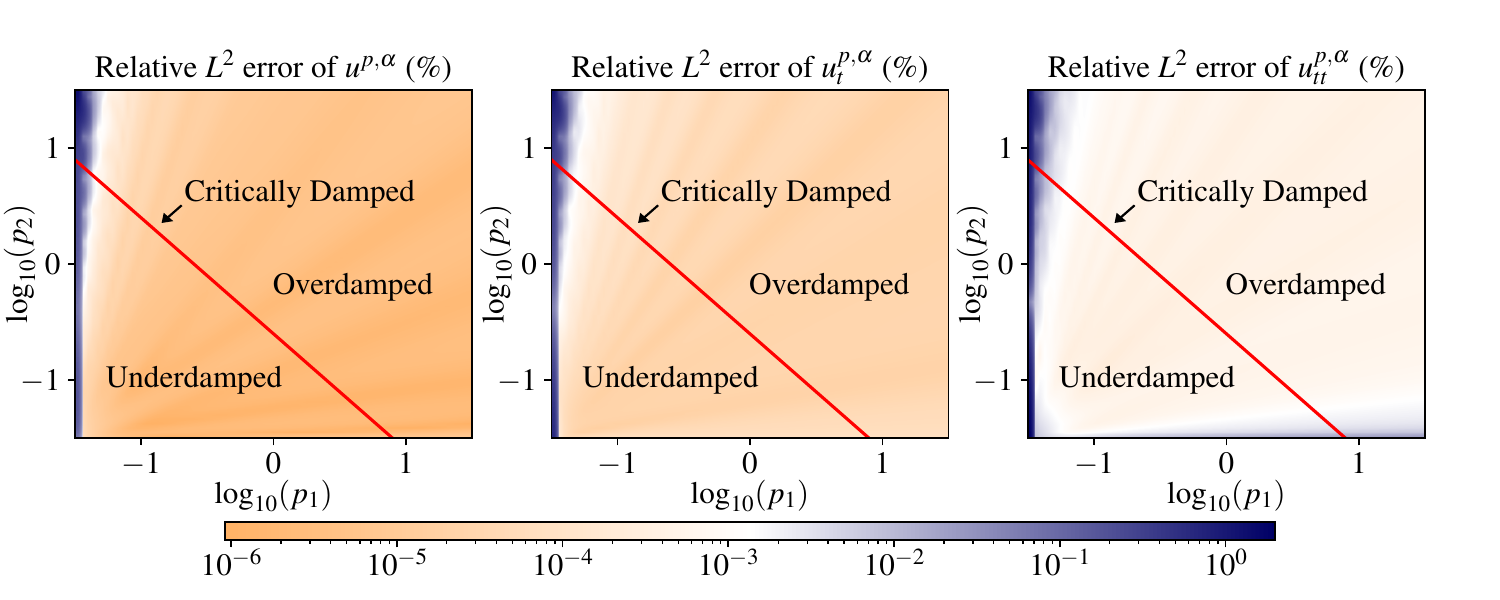}
	\caption{After training.}\label{fig_PINN_L2_Error_b}
	\end{center}
    \end{subfigure}

	\end{center}
			\vspace{-0.4cm}
	\caption{Relative {$L^2$-errors} (in $\%$) of the LS-Net solutions and their first and second derivatives for the damped harmonic oscillator problem. We display the errors after solving the corresponding LS problem, (a) before and (b) after training the NN. {Note the difference in the ranges of the color bars.}}\label{fig_PINN_L2_Error}
\end{figure}

Figure \ref{fig_PINN_L2_Error} presents the relative {$L^2$-errors} (in $\%$) of the LS-Net solutions $u^{p, \alpha}$ and their first and second derivatives by solving the corresponding LS problem before and after training the NN. {To present the results of this test, the relative errors are computed for a $60 \times 60$ grid of parameters $p_1$ and $p_2$, logarithmically distributed within the range $[10^{-1.5}, 10^{1.5}]$.} Figure \ref{fig_PINN_L2_Error_a}  demonstrates that, before training the NN, the error range varies approximately from $10^{-2} \%$ to $10^{3}\%$. However, the majority of the data exhibit relative errors around  $100\%$ for $u^{p, \alpha}$ and its first and second derivatives. Additionally, we observe that the errors increase to nearly $1000\%$ as $p_1$ approaches zero. A small value of $p_1$ corresponds to a convected-dominated diffusion-type ODE, leading to singular behavior. Examination of Figure \ref{fig_PINN_L2_Error_b} reveals a substantial improvement in the accuracy across the entire parameter range. Indeed, after training, we have typical relative errors of orders $10^{-6}\%$, $10^{-5}\%$, and $10^{-3}\%$, respectively, for $u^{p, \alpha}$ and its first and second derivatives. Additionally, we observe that, after training, the relative {$L^2$-errors} for values of $p_1$ close to zero are on the order of $1\%$, which is satisfactory.

\begin{figure}[H]
	\begin{center}
    \includegraphics[width=1\textwidth]{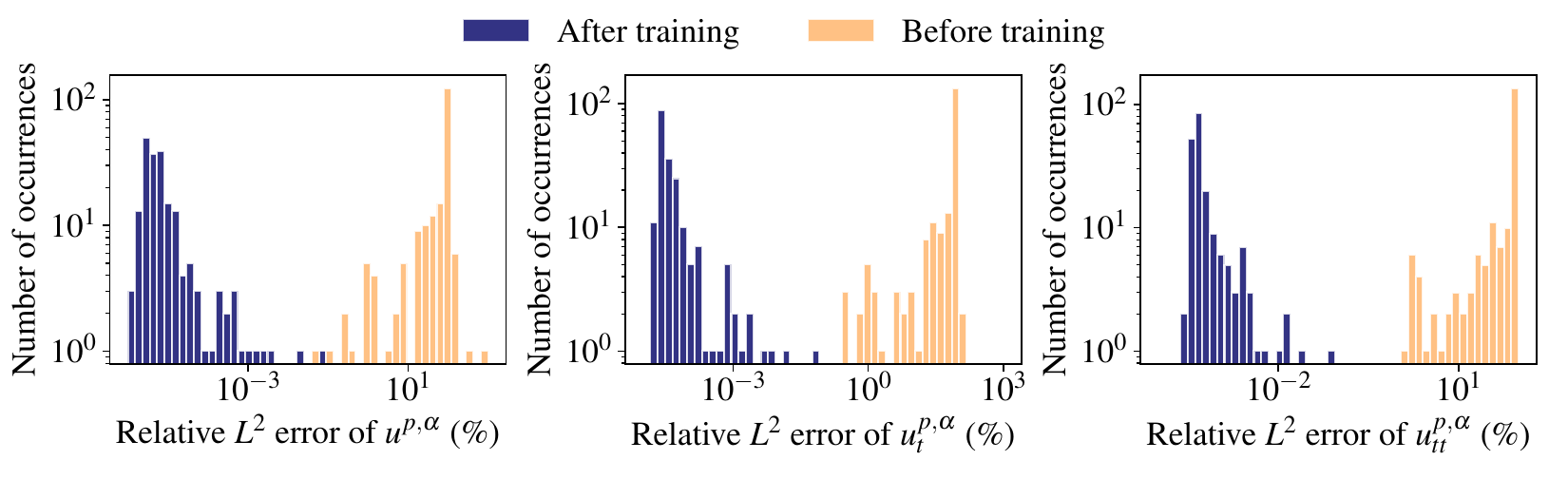}
    \end{center}
    \vspace{-0.7cm}
	\caption{Distribution of the relative {$L^2$-errors} (in $\%$) for the LS-Net solution and their first and second derivatives of the damped harmonic oscillator problem. Results include solving the LS problem before and after training the NN.}\label{fig_PINN_Distribution}
\end{figure}

\begin{figure}[H]
	\begin{center}
    \includegraphics[width=1\textwidth]{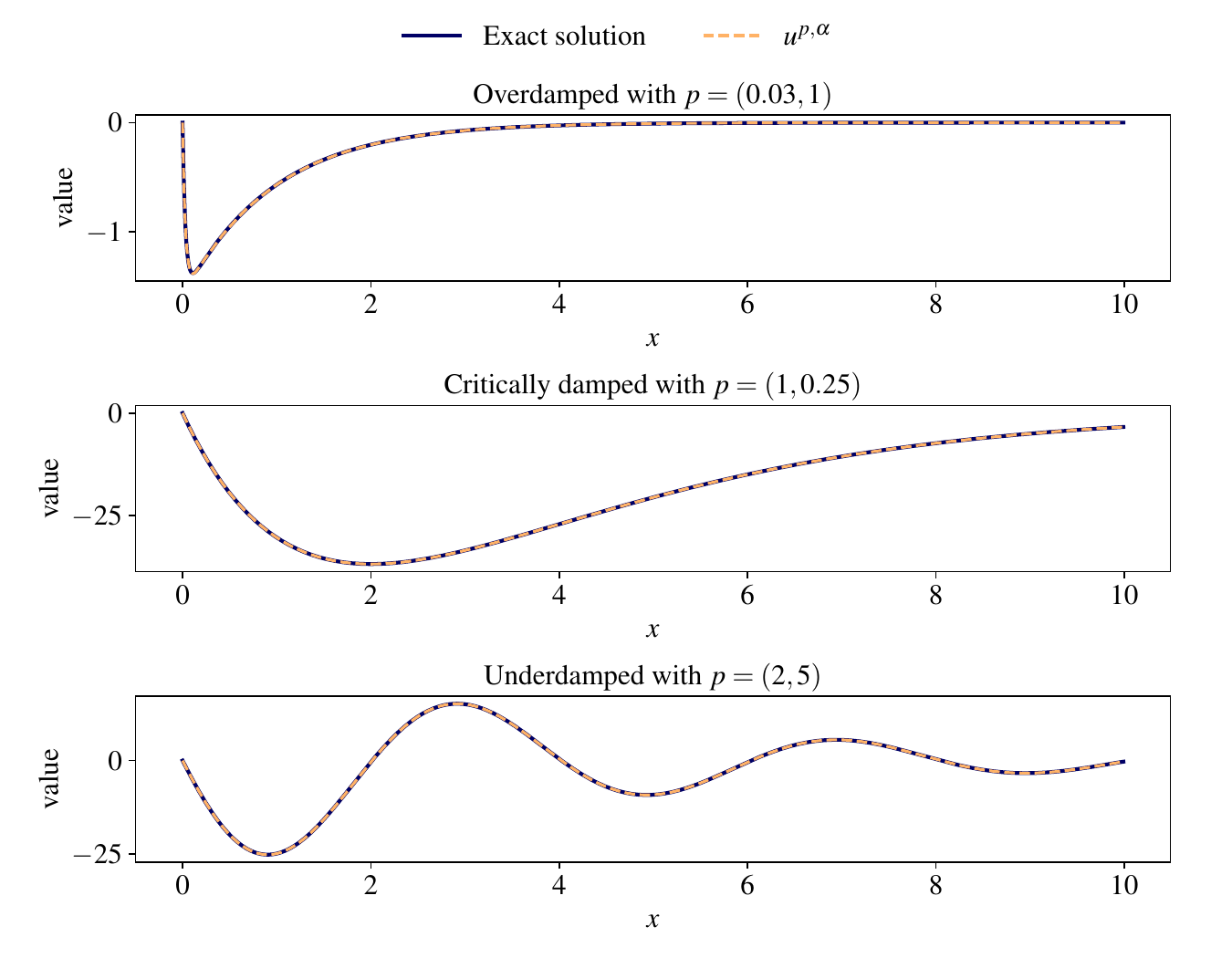}
    \end{center}
    \vspace{-0.7cm}
	\caption{The exact and the LS-Net solutions of the damped harmonic oscillator problem, illustrate the three damping behaviors: {overdamped, critically damped, and underdamped.}}\label{fig_PINN_exact_LS}
\end{figure}

To better illustrate the distribution of the relative  {$L^2$-errors} both before and after training the NN, we considered a parameter set of $10^3$ values sampled log-uniformly and uniformly, respectively, for $p_1$ and $p_2$, over the interval $[10^{-1.5}, 10^{1.5}]$.  Figure \ref{fig_PINN_Distribution} illustrates the results for $u^{p, \alpha}$ and its first and second derivatives. As seen, before training we have the typical relative {$L^2$-error} distribution of order $50\%$ for $u^{p, \alpha}$ and its first and second derivatives.  In contrast, after training, the typical relative {$L^2$-errors} are reduced to approximately  $10^{-5}\%$, $10^{-4}\%$, and $10^{-3}\%$ for $u^{p, \alpha}$ and its first and second derivatives, respectively. Figure \eqref{fig_PINN_exact_LS} presents the exact and LS-Net solutions of problem \eqref{eq_damped_problem_2p} for three different parameter sets, illustrating the three possible damping behaviors.

\begin{figure}[h]
    \centering
        \includegraphics[width=0.7\textwidth]{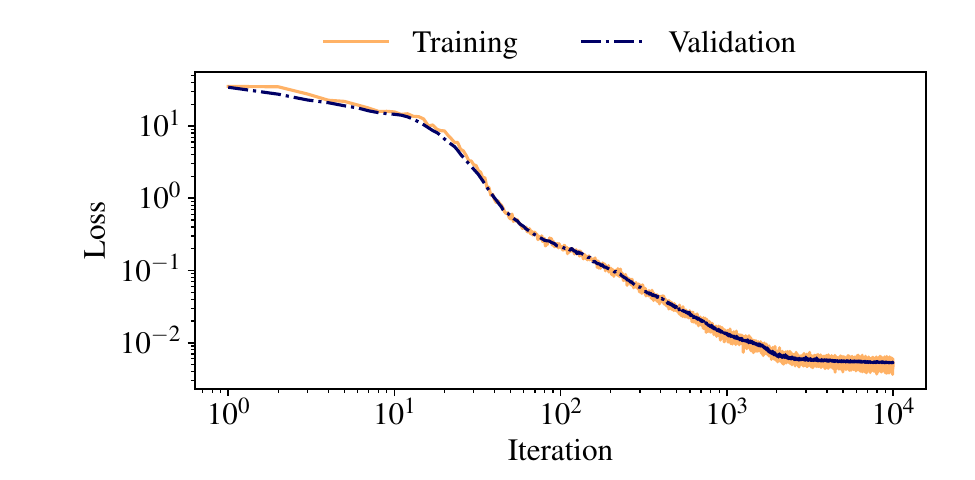}
        \vspace{-0.1cm}
        \caption{Loss evolution of the LS-Net method on the training and validation datasets for the damped harmonic oscillator problem.}\label{fig_PINN_Damping_Loss}
    
\end{figure}

{Finally, Figure \ref{fig_PINN_Damping_Loss} displays the loss evaluation of the LS-Net method on the training and validation datasets over $10^{4}$ training iterations. Recall that the training loss is calculated using 500 parameters randomly sampled at each iteration, while the validation loss is assessed on a fixed set of 500 parameters, preselected before training begins.} The validation loss values are in the mid-range of the training loss values, showing no noticeable overfitting present.


\subsection{One-dimensional Helmholtz equation}\label{One-dimensional Helmholtz equation}
Let $\Omega = (0, 1)$. We consider the one-dimensional parametric Helmholtz equation governed by the following weak form with homogeneous impedance boundary condition:
find $\mathfrak{u}^p\in H^1(\Omega)$, such that for all $v\in H^1(\Omega)$, the following holds
\begin{equation}
    \displaystyle \int_\Omega  \left( p_1 \mathfrak{u}^p_x v_x  - p_2^2 \mathfrak{u}^pv +  fv\right)  \ dx - i \sqrt{p_1} p_2 \int _{\partial \Omega}  \mathfrak{u}^pv \ dx =0,
\end{equation}
where $i$ denotes the imaginary unit. In applications, $p_2 = \omega /c$ such that $\omega$ is the angular frequency and $c$ is the wave speed. Here, we restrict ourselves to the case when $c$ is a positive constant. 

\noindent Additionally,  we define $f$  as 
\begin{equation}
f(x)=\left\{\begin{array}{ll}1, & 0\leq x\leq 0.5,\\
0, & 0.5 \leq x\leq 1,\end{array}\right.
\end{equation} 
and consider the discontinuous parameter $p:(0,1)\to\mathbb{R}$ as
\begin{equation}\label{eq_sigma}
p_1(x)=\left\{\begin{array}{ll}p_{11}, & 0\leq x\leq 0.5,\\
p_{12}, & 0.5 \leq x\leq 1.\end{array}\right.
\end{equation}
The exact solution of the considered problem is given by
\begin{equation}\label{eq_ue_Complex-valued_Helmholtz}
 \mathfrak{u}^p(x)=\left\{\begin{array}{ll}a \exp\left(-i\displaystyle \frac{p_2 x}{\sqrt{p_{11}}}\right), &\displaystyle 0\leq x\leq 0.5,\\[4mm]
b \exp\left(i \displaystyle\frac{p_2 x}{\sqrt{p_{12}}}\right) + c \exp\left(-i \displaystyle \frac{p_2 x}{\sqrt{p_{12}}}\right) + \displaystyle\frac{1}{p_2^2}, &\displaystyle 0.5 \leq x\leq 1,\end{array}\right.
\end{equation}
where 
\[
\begin{array}{rl}
    a =& \hspace{-2mm} \displaystyle \frac{\Big(-\exp\left(\displaystyle\frac{\frac{i}{2}p_2(\sqrt{p_{11}} + \sqrt{p_{12}})}{\sqrt{p_{11} p_{12}}}\right) + \exp \left(\displaystyle\frac{\frac{i}{2} p_2}{\sqrt{p_{11}}}\right) \Big) p_{12} \sqrt{p_{11}}}{p_2 ^2 \left(\sqrt{p_{12}}p_{11} + \sqrt{p_{11}}p_{12}\right)}, \\[6mm]
    
    b = &\hspace{-2mm}\displaystyle \frac{-2 p_{11} \sqrt{p_{12}} \exp \left(\displaystyle\frac{-\frac{i}{2} p_2}{\sqrt{p_{12}}} - \sqrt{p_{11}}p_{12} + \sqrt{p_{12}}p_{11}\right)}{2p_2 ^2 \left(\sqrt{p_{12}}p_{11} + \sqrt{p_{11}}p_{12}\right)}, \\[6mm]
    
    c =& \hspace{-2mm}\displaystyle - \frac{\exp \left(\displaystyle\frac{ip_2}{\sqrt{p_{12}}}\right)}{2 p_2 ^2}. 
\end{array}
\]

{We employ an NN of the form described in \eqref{eq_regularity-conforming} with $K=3$, $N_1 = N_2 = 5$, and $N_3 =30$.}
Given that this problem involves a piecewise-constant parameter $p_1$ with a discontinuity at $x = 0.5$, the derivative of $\mathfrak{u}^p$ at this point will generally be discontinuous. This discontinuity poses challenges when attempting to approximate $\mathfrak{u}^p$ by a smooth function. To this end, consider the function $\eta(x) = |x-0.5|$, which is a Lipschitz function with a discontinuity in its derivative at $x = 0.5$.  Then,  function $\psi$ is defined as 
\begin{equation}\label{eq_complex_phi}
{\psi = (\underbrace{1, \ldots, 1}_{15}, \underbrace{\eta, \ldots, \eta}_{15})^T.}
\end{equation}
This structure enables us to enforce the interface condition, ensuring that the left and right limits of $p_1 u^{p, \alpha}_x$ are equal.
Due to the homogeneous impedance boundary condition, the NN structure under consideration has $\varphi = 1$. 
{A discretized loss function based on the DFR approach, using 400 test functions and 1,000 quadrature points, was employed to train $\mathbf{u}^\alpha$.} The training set consists of 500 samples for each parameter, where $p_{11}$ and $p_{12}$  are drawn from a log-uniform distribution over the interval  $(0.05, 100]$ and the parameter $p_2$ is sampled from a uniform distribution over $(0, 10]$. {The validation set is constructed using the same sampling approach, consisting of 500 samples for each parameter, fixed prior to the start of training. Additionally, we use 1000 quadrature points to compute the training loss, while the validation loss is evaluated using 1005 quadrature points. } The network is trained for $5 \times 10 ^{3}$ iterations using the Adam optimizer, with a learning rate that exponentially decreased from $10^{-1}$ to $10^{-4}$. 

\begin{remark}\label{DFR test function 1D}
The eigenvectors of the Laplacian, subject to appropriate homogeneous boundary conditions, are considered as test functions in the DFR method. In one-dimensional problems, the eigenvectors are known as cosines/sines. Herein, we have $v_m = a_m \cos((m-1)\pi x)$, where for $m=1, \ldots, M$ the constants $a_m$ are such that $||v_m||_{H^1}=1$. 
\end{remark}

\begin{figure}[H]
    \centering
        \includegraphics[width=0.7\textwidth]{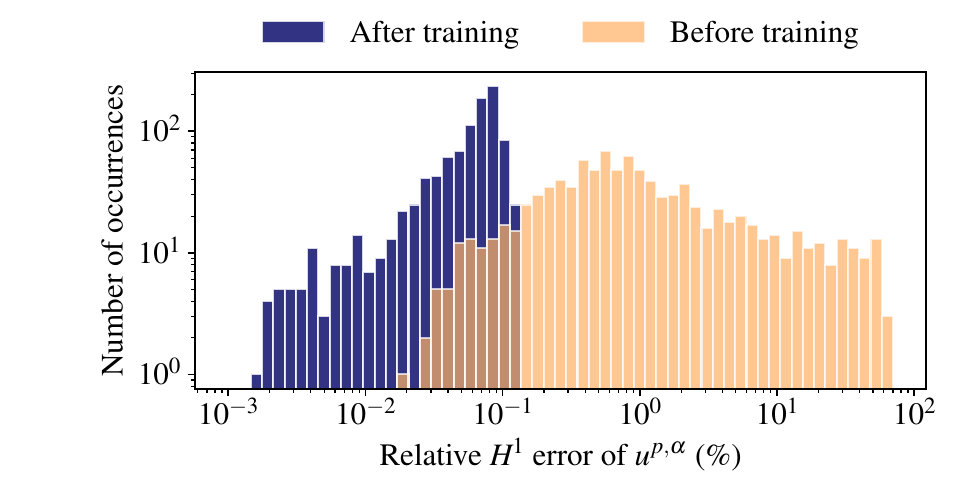} 
        \caption{Distribution of the relative $H^1$-errors (in $\%$) for the LS-Net solution of the one-dimensional Helmholtz equation. The results are obtained by solving the corresponding LS problem before and after training the NN.}\label{fig_Distribution_HH}
    
\end{figure}

\begin{figure}[H]
    \centering
    \begin{subfigure}[b]{0.9\textwidth}
        \centering
        \includegraphics[width=\textwidth]{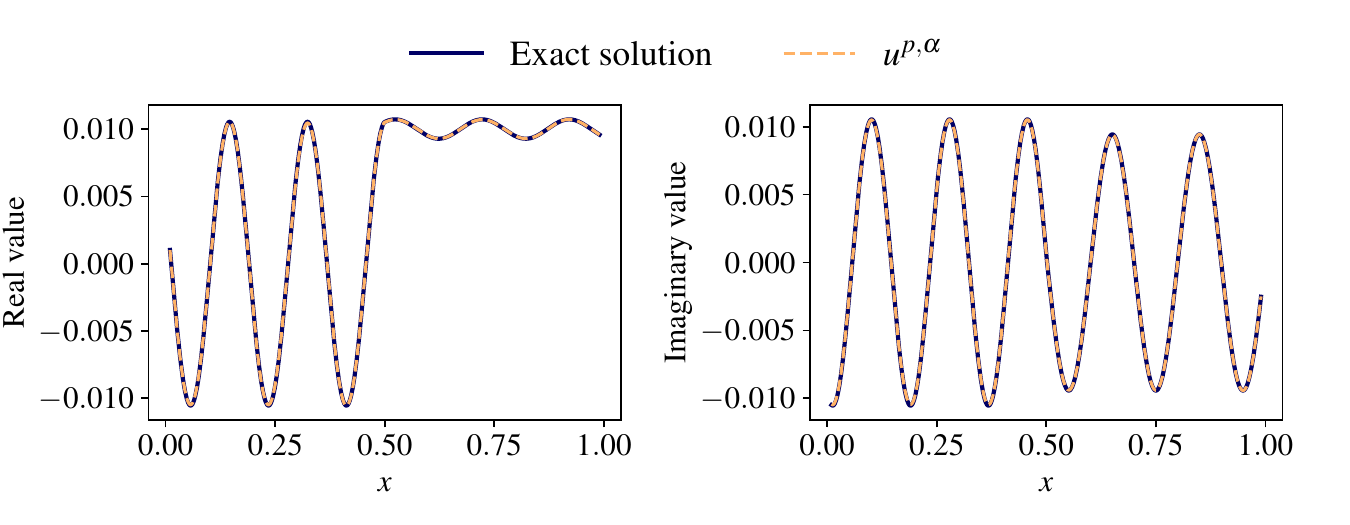}
        \caption{$p = (0.08, 0.1, 10)$}
        \label{fig:high_freq}
    \end{subfigure}
    
    \begin{subfigure}[b]{0.9\textwidth}
        \centering
        \includegraphics[width=\textwidth]{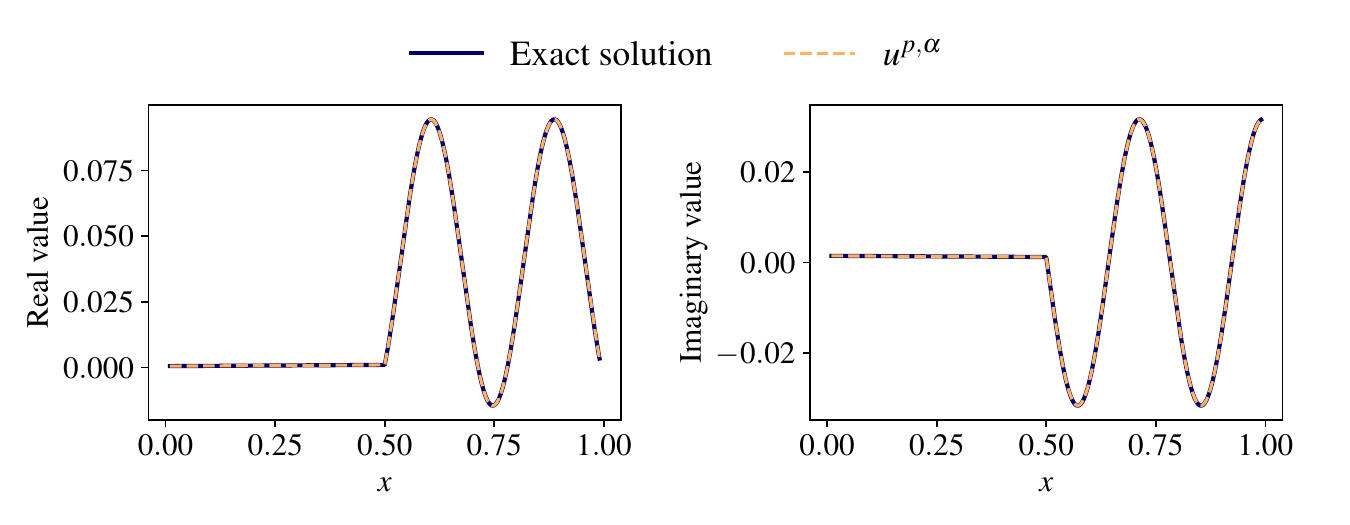}
        \caption{$p = (50, 0.05, 5)$}
        \label{fig:low_freq}
    \end{subfigure}
    
    \caption{The exact and the LS-Net solution of the one-dimensional Helmholtz equation for two different values of parameters.}
    \label{fig_high_low_freq}
\end{figure}

\begin{figure}[H]
    \centering
        \includegraphics[width=0.7\textwidth]{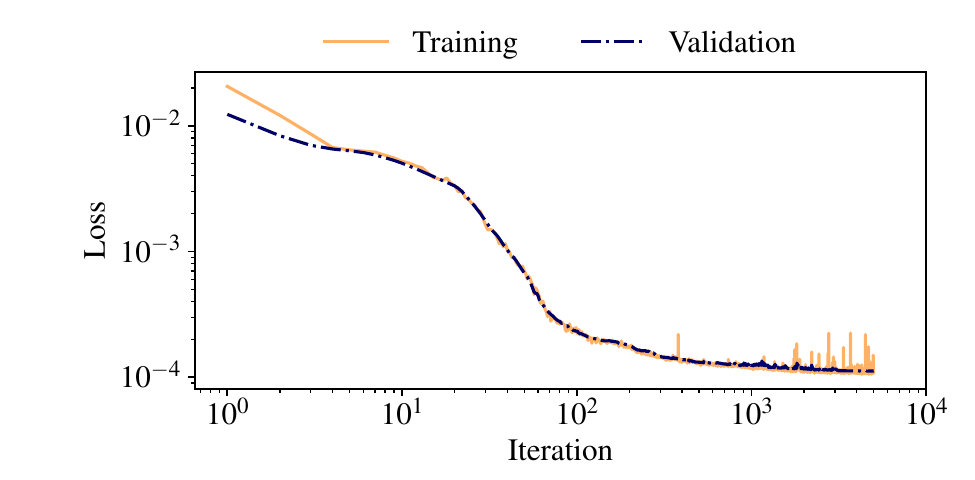} 
        \caption{Loss evolution of the LS-Net method on the training and validation datasets for the one-dimensional Helmholtz equation.}\label{fig_HH_Loss}
    
\end{figure}

Figure \ref{fig_Distribution_HH} illustrates the relative $H^1$-errors (in $\%$) of the LS-Net solutions by solving the corresponding LS problems before and after training the NN. As shown, this indicates the effectiveness of the training process in significantly reducing the relative $H^1$-errors. Specifically, the minimum and maximum relative errors decrease from approximately $0.8\%$ and $70\%$ to $0.009\%$ and $0.1\%$, respectively. Figure \ref{fig_high_low_freq} also shows the exact \eqref{eq_ue_Complex-valued_Helmholtz} and LS-Net solutions for two sets of parameters $p=(0.08, 0.1, 10)$ and $p=(50, 0.05, 5)$. The results demonstrate that the LS-Net method effectively approximates high-frequency solutions and accurately captures the singularities in the exact solutions. Furthermore, Figure \ref{fig_HH_Loss} displays the loss evaluations for the LS-Net method on both the training and validation datasets throughout the training process. We see a strong alignment between the training and validation losses, indicating that no overfitting has occurred.

\subsection{Two-dimensional transmission problem}\label{Two-dimensional transmission problem}

\begin{figure}[h]
\begin{center}
\scalebox{0.95}{\definecolor{My_blue}{RGB}{0,0,102}
\definecolor{My_Orange}{RGB}{255, 178, 102}
\begin{tikzpicture}[
punkt/.style={
   rectangle,
   rounded corners,
   draw=black, thick,
   text width=2em,
   minimum height=1em,
   text centered}
]
    \filldraw[fill=My_Orange!40!white, thick] (-2.2,-2.2) rectangle (2.2,2.2);
    
    \node at (2.7, 0) [right] {$x$};
    \node at (0, 2.7) [above] {$y$};

    \draw[->] (-2.7,0) -- (2.7,0) node[anchor=north west] {};
    \draw[->] (0,-2.7) -- (0,2.7) node[anchor=south east] {};

    \draw[fill=My_blue!45!white, thick] (-1.1, 1.1) circle (0.55);
    \draw[fill=My_blue!45!white, thick] (1.1, 1.1) circle (0.55);
    \draw[fill=My_blue!45!white, thick] (-1.1, -1.1) circle (0.55);
    \draw[fill=My_blue!45!white, thick] (1.1, -1.1) circle (0.55);

    \draw[thin, dashed] (-1.1,-1.1) rectangle (1.1,1.1);


    \filldraw (-2.2, 0) circle (1pt) node[anchor=north east] {$-1$};
    \filldraw (2.2, 0) circle (1pt) node[anchor=north west] {$1$};
    \filldraw (0, -2.2) circle (1pt) node[anchor=north west] {$-1$};
    \filldraw (0, 2.2) circle (1pt) node[anchor=south west] {$1$};

    \filldraw (-1.1, 0) circle (1pt) node[anchor=north east] {$-\frac{1}{2}$};
    \filldraw (1.1, 0) circle (1pt) node[anchor=north west] {$\frac{1}{2}$};
    \filldraw (0, -1.1) circle (1pt) node[anchor=north west] {$-\frac{1}{2}$};
    \filldraw (0, 1.1) circle (1pt) node[anchor=south west] {$\frac{1}{2}$};

    \node at (1.5, 1.5) [right] {$\Omega_1$};
\node at (-2.2, 1.5) [right] {$\Omega_2$};
\node at (1.5, -1.5) [right] {$\Omega_3$};
\node at (-2.2, -1.5) [right] {$\Omega_4$};
\node at (-0.3, -3) [right] {$\Omega$};

\end{tikzpicture}}
\end{center}	
\vspace{-0.5cm}
\caption{The geometry of the two-dimensional domain $\Omega$ for transmission problem.}\label{fig_geometry_2D}
\end{figure}
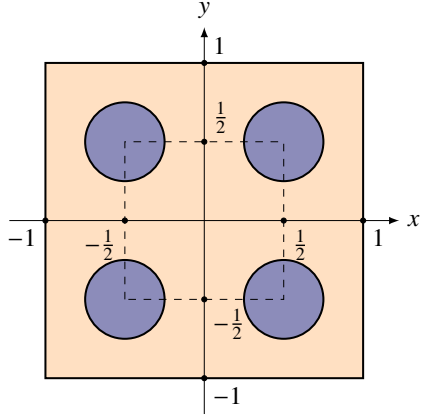

Let $\Omega = (-1,1)^2$. Assume that $\Omega_i$ are circles in the domain $\Omega$, each with a radius of $\frac{1}{4}$ and centered at the points $(x_i,y_i)=\left(\pm\frac{1}{2},\pm\frac{1}{2}\right)$ for $1\leq i\leq 4$. The geometric configuration is illustrated in Figure \ref{fig_geometry_2D}.
We consider the transmission problem governed by the following weak-form Poisson's equation with inhomogeneous Dirichlet boundary data: find $\mathfrak{u}^p\in  \cos\left(\frac{\pi}{2}x\right) +H^1_0(\Omega)$, such that 
\begin{equation}
\int_\Omega p\nabla \mathfrak{u}^p\cdot\nabla v = 0,\qquad \forall v\in H^1_0(\Omega),
\end{equation}

\noindent where the parameter $p$ is defined as piecewise constant in $\Omega$ by
\begin{equation}
p(x,y)=\left\{\begin{array}{l l}p_i, & (x,y)\in \Omega_i,\\ 1, & \text{otherwise}, \end{array}\right.\qquad p_i \in \mathbb{P}=[1,10].
\end{equation} 
We use the DFR method equipping the test space by $\Vert\cdot\Vert_{H^1_0(\Omega)}:=\Vert\nabla\cdot\Vert_{L^2(\Omega)}$ and choosing corresponding orthonormal sinusoidal functions for the discretization. {For the trial space, we consider an NN with $K = 4$ and $N_k = 75$ for $k=1, \ldots, 4$.}
The boundary conditions and regularity are enforced using the cut-off functions $\varphi(x, y) = (x^2-1)(y^2-1)$ and
\begin{equation}\label{eq_2D_phi}
\psi = (\underbrace{1, \ldots, 1}_{15}, \underbrace{\eta_1, \ldots, \eta_1}_{15}, \underbrace{\eta_2, \ldots, \eta_2}_{15}, \underbrace{\eta_3, \ldots, \eta_3}_{15}, \underbrace{\eta_4, \ldots, \eta_4}_{15})^T,
\end{equation}
where 
$$\eta_i(x, y) = \text{LeLU}\left(\frac{1}{4} - ||(x_i,y_i)-(x,y)||^2_2\right), \qquad 1\leq i \leq 4,$$
and the LeLU (Linear Exponential Linear Unit) function is given by
\begin{equation}
\text{LeLU}(t)=\left\{\begin{array}{l l}t, & t\geq 0,\\ 2t, & \text{otherwise}. \end{array}\right.
\end{equation}
We highlight that $u^{p, \alpha}$ is designed to approximate the solution of the associated homogeneous problem in $H^1_0(\Omega)$. Then, the final solution is obtained by adding the corresponding lift, i.e., $u^{p, \alpha} + \cos\left(\frac{\pi}{2}x\right)$.

We train for $10^4$ iterations using the Adam optimizer with an initial learning rate of $10^{-2}$. We consider $75\times 75$ sinusoidal basis functions for the test space discretization, $300\times 300$ fixed quadrature points for integration, and a batch size of $32$ parameter values that are sampled following uniform random variables iteratively in discretized loss function \eqref{eq_Loss-D} during training. We consider two settings for computing the loss on the validation dataset: (i)  integration validation, which uses the same number of test basis functions as during training  ($75\times 75$) but with a larger number of quadrature points ($400\times 400$), and (ii) truncation validation, which employs a greater number of both test basis functions ($100\times 100$) and quadrature points ($400\times 400$). In both cases, we consider a fixed batch size of $500$ parameter values.

\begin{remark}
    In tensor product domains, the tensor products of one-dimensional eigenvectors yield an orthogonal basis of eigenvectors for the Laplacian. Specifically, we have $v_m = a_m \sin(m_1\pi x)\sin(m_2\pi x)$, where for $m = (m_1, m_2)$ and $m_1, m_2=1, \ldots, M$ the constants $a_m$ are such that $||v_m||_{H^1}=1$. 
\end{remark}
\noindent As the exact parametric solution $\mathfrak{u}^p$ presents a non-elementary, closed formula that becomes challenging to treat in practice, we use the following lower and upper bounds of the relative error
\begin{equation}\label{eq_lower and upper}
\frac{1/\vartheta^p  {\sqrt{\mathcal{L}^p(\alpha)}}}{\Vert u^{\alpha,\mathbf{c}}\Vert_{H^1_0(\Omega)} + 1/\vartheta^p \sqrt{\mathcal{L}^p(\alpha)}} \leq \frac{\Vert u^{\alpha,\mathbf{c}}-\mathfrak{u}^p\Vert_{H^1_0(\Omega)}}{\Vert \mathfrak{u}^p\Vert_{H^1_0(\Omega)}} \leq \frac{{\sqrt{\mathcal{L}^p(\alpha)}}}{\Vert u^{\alpha,\mathbf{c}}\Vert_{H^1_0(\Omega)} - \sqrt{\mathcal{L}^p(\alpha)}},
\end{equation} 
where $\vartheta^p := \max\{p_1,p_2,p_3,p_4\}$ is a per-parameter estimate for the continuity (sup-sup) constant---the estimate for the coercivity (inf-sup) constant is one. In practice, we approximate $\Vert u^{\alpha,\mathbf{c}}\Vert_{H^1_0(\Omega)}$ and $\mathcal{L}^p(\alpha)$ employing the truncation validation setting specified above.

\begin{figure}[H]
    \centering
        \includegraphics[width=0.97\textwidth]{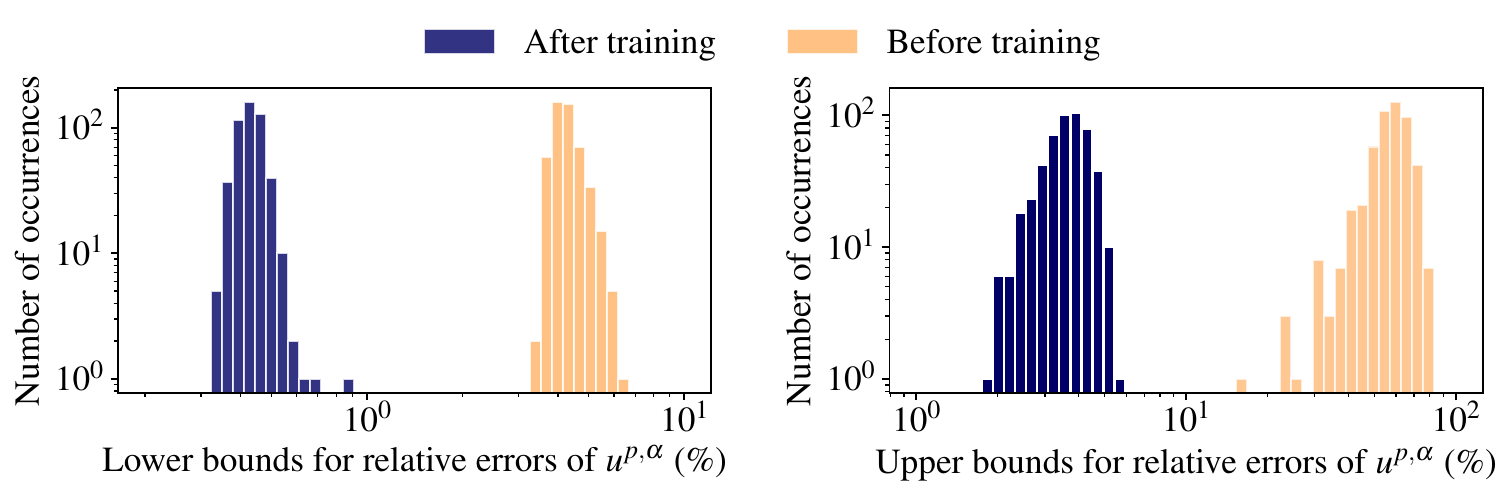} 
        \caption{Distribution of the upper and lower bounds \eqref{eq_lower and upper} for the relative {$L^2$-errors} (in $\%$)  of the LS-Net solution for the two-dimensional transmission problem. Results are obtained by solving the corresponding LS problem before and after training the NN.}\label{fig_Distribution_transmission}
    
\end{figure}

\begin{figure}[H]
    \centering
        \includegraphics[width=0.85\textwidth]{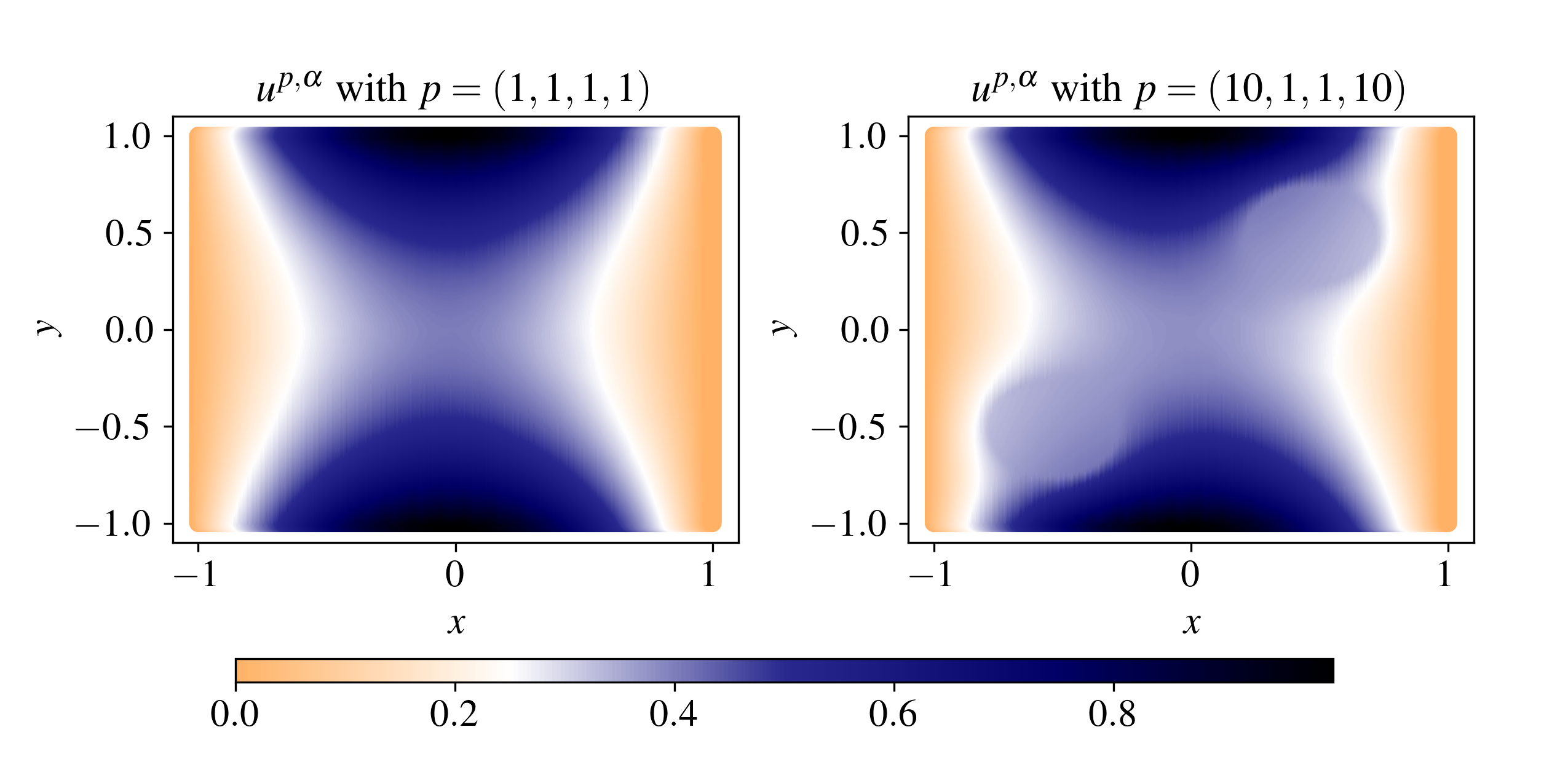} 
        \vspace{-0.4cm}
        \caption{LS-Net solution $u^{p, \alpha}$ of the two-dimensional transmission problem for two different sets of parameters.}\label{fig_2D_LS-Net_solution}
    
\end{figure}

\begin{figure}[h]
    \centering
        \includegraphics[width=0.65\textwidth]{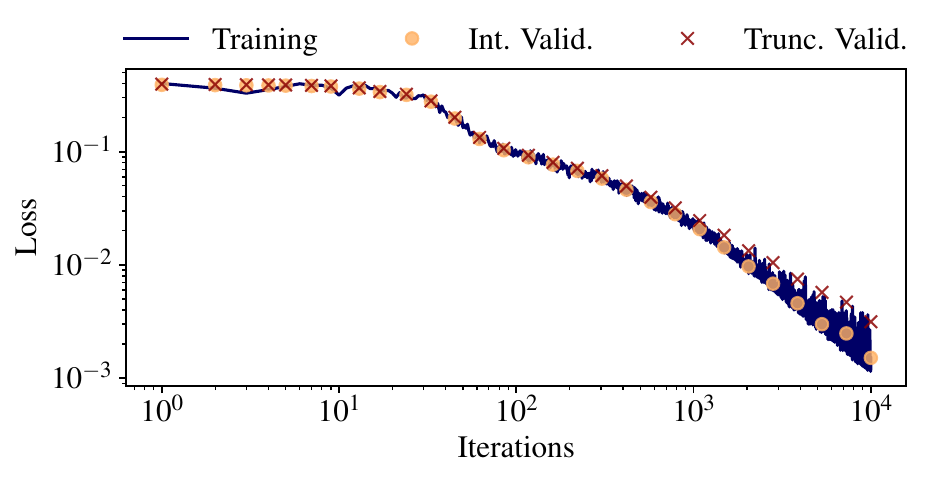} 
        \caption{Loss evolution of the LS-Net method for the two-dimensional transmission problem on the training dataset and two validation levels.}\label{fig_Loss_transmission}
    
\end{figure}

Figure \ref{fig_Distribution_transmission} illustrates the lower and upper bounds for the relative $H^1$-errors (in $\%$) of the LS-Net solution by solving the corresponding LS problems before and after training the NN. For each parameter, we considered a set of $500$ values sampled uniformly over the interval $[1, 10]$.
Prior to training, the majority of the lower bounds fall within the range of $[2\%, 5\%]$, while the upper bounds are clustered between $[30\%, 80\%]$. After training, both the upper and lower bounds significantly decrease, with the majority now concentrated in narrower intervals  $[0.3\%, 0.7 \%]$ and $[1\%, 5\%]$, respectively, indicating improved accuracy. Figure \ref{fig_2D_LS-Net_solution} displays the LS-Net solutions for two parameter sets $p=(1, 1, 1, 1)$ and $p=(10, 1, 1, 10)$.  Moreover, Figure \ref{fig_Loss_transmission} illustrates the loss evolution of the training dataset, along with computing the integration and truncation validations throughout the training process. The figure shows a close alignment between the training and validation loss values, indicating a well-converged and robust training process.

\vspace{-0.7cm}
\vspace{-0.2cm}
\section{Conclusions}\label{Conclusion}
In this work, we introduced the LS-Net method as a novel approach for solving linear parametric PDEs. This method employs a single deep NN and an LS solver to generate parametric solutions in a separated representation form. Specifically, the deep NN generates a vectorial function, where its components serve as basis functions for the solution space. The LS solver then computes the coefficients of the basis functions for each parameter distribution. The LS-Net method can be interpreted as a variant of the PGD method, effectively addressing the challenges encountered in classical approaches by utilizing the strengths of NNs in approximating complex functions. To pose the LS problem, the LS-Net method requires a quadratic loss function. In this work, two loss functions based on the DFR and PINN approaches have been used. We provided some theoretical results in style of a Universal Approximation Theorem, stating that the LS-Net method can achieve the desired accuracy by utilizing sufficiently large NNs in combination with an LS solver. 

We also validated the obtained theoretical results by solving three benchmark problems. First, the damped harmonic oscillator problem was solved, where LS-Net effectively approximated various damping behaviors and achieved a significant reduction in errors after training the NN. Next, we applied the LS-Net method to a one-dimensional Helmholtz problem with impedance boundary conditions. The LS-Net performed well again, particularly in managing solution discontinuities through a regularity-conforming NN. Finally, a two-dimensional transmission problem has been solved. In this example, due to the absence of an analytical solution, we employed the upper and lower bounds for the relative $H^1$-errors of the LS-Net solution.

Although the current LS-Net method is designed for linear PDEs, the concept may be extended to non-linear PDEs by considering the linearization of the PDE analogously to a Newton method, although providing well-chosen initial guesses, which may generally be parameter-dependent, may prove challenging. Moreover, the current LS-Net method necessitates solving an LS system for each new parameter, which can lead to computational overhead, especially as the number of integration points or basis functions grows. This issue is partially addressed by performing a single high-precision integration and constructing the LS system parametrically after training the network. This approach eliminates the need for costly integrations at each iteration, making the process significantly more efficient. Our future work will aim to further optimize this process to reduce computational complexity and enhance overall efficiency.

{\section*{Acknowledgments}}\label{Acknowledgments}
This work has received funding from the following Research Projects/Grants: European Union’s Horizon Europe research and innovation programme under the Marie Sklodowska-Curie grant agreement No 101119556 (IN-DEEP). TED2021-132783B-I00 funded by MICIU/AEI /10.13039/501100011033 and by FEDER, EU; PID2023-146678OB-I00 funded by MICIU/AEI /10.13039/501100011033 and by the European Union NextGenerationEU/ PRTR; “BCAM Severo Ochoa” accreditation of excellence CEX2021-001142-S funded by MICIU/AEI/10.13039/5011000; Misiones Project IA4TES MIA.2021.M04.0008 funded by the Ministry of Economic Affairs and Digital Transformation; Basque Government through the BERC 2022-2025 program; BEREZ-IA (KK-2023/00012) and RUL-ET (KK-2024/00086), funded by the Basque Government through ELKARTEK; Consolidated Research Group MATHMODE (IT1456-22) given by the Department of Education of the Basque Government; BCAM-IKUR-UPV/EHU, funded by the Basque Government IKUR Strategy and by the European Union NextGenerationEU/PRTR.

\vspace{0.5cm}
\section*{Appendix: Universal Approximation Theorem for the LS-Net method}\label{Appendix}

{For simplicity, we denote the Sobolev space  $H^k(\Omega)$ by $H$ for $k\in \mathbb{N}$. We also assume the existence of the operator, 
\begin{equation}\label{eq_operatpr_A}
A:=\int_\mathbb{P} \mathfrak{u}^p\otimes \mathfrak{u}^p\,d\mu(p),
\end{equation}
understood as a Bochner integral in the space $L^1(\mathbb{P},B(H);\mu)$, where $B(H)$ is the set of bounded operators from $H\to H$ with the operator norm.  The Bochner integral generalizes the Lebesgue integral to functions that take values in Banach spaces such as $B(H)$ \cite{hille1996functional}. The operator $A$ acts on elements $v\in H$ by 
\begin{equation}
Av = \left(\int_\mathbb{P} \mathfrak{u}^p\otimes \mathfrak{u}^p\,d\mu(p)\right)v=\int_\mathbb{P} \langle \mathfrak{u}^p,v\rangle_H \mathfrak{u}^p\,d\mu(p), 
\end{equation}
where $\langle \cdot, \cdot \rangle_H$ denotes the standard inner product of the Sobolev space  $H$. Therefore, $A(v)$ belongs to  $H$, as  $\langle \mathfrak{u}^p,v\rangle_H$ produces a scalar for each $p\in \mathbb{P}$, and multiplying this scalar by $\mathfrak{u}^p \in H(\Omega)$  results in an element of $H$.  Moreover, the Bochner integral over $\mathbb{P}$ with respect to $\mu$ ensures that the result remains in $H$, as $H$ is a closed vector space under such integrals. On the other hand, using the triangle inequality for the norm and the Cauchy-Schwarz inequality in $H$, we have
\begin{equation}
\left\Vert A(v) \right\Vert _H = \left\Vert \int_\mathbb{P} \langle \mathfrak{u}^p,v\rangle_H \mathfrak{u}^p\,d\mu(p) \right\Vert _H \leq \left\Vert v \right\Vert_H  \int_\mathbb{P} \left\Vert \mathfrak{u}^p\right\Vert ^2_H \,d\mu(p)< \infty.
\end{equation}
This inequality confirms that $A(v)$ has a finite $H$-norm. Therefore, the Bochner integral ensures that the operator $A$ is well-defined as an integral of bounded operators acting on $H$.}


\begin{proposition}\label{prop.AppNNSubspace}
Let $H_N\subset H$ be a finite-dimensional subspace of $H$, with dimension $N$. Then, for any basis $\{u_n\}_{n=1}^N$ of $H_N$ and $\epsilon>0$, there exists a fully-connected feed-forward NN $\mathbf{u}^\alpha:\mathbb{R}^{N_0}\to\mathbb{R}^N$ such that $\displaystyle\sum\limits_{n=1}^N||u^\alpha_n-u_n||_{H}<\epsilon$. 
\end{proposition}
\begin{proof}
This is a consequence of the Universal Approximation Theorem proved in \cite{hornik1990universal}. 
\end{proof}

\begin{proposition}\label{prop.AppSubspaceL2}
For every $\epsilon>0$, there exists a finite-dimensional subspace $H_N$ of $H$ such that 
\begin{equation}
\int_\mathbb{P} \min\limits_{u\in H_N}||u-\mathfrak{u}^p||_H^2\,d\mu(p)<\epsilon. 
\end{equation}
\end{proposition}

\begin{proof}
First, we consider $H_N$ to be an arbitrary, finite dimensional subspace of $H$, with orthonormal basis $\{u_n\}_{n=1}^N$. 
We denote $\mathscr{P} _{_{H_N}}:H\to H_N$ the orthogonal projection operator from $H$ onto the subspace $ H_N $. We may then express the integral in terms of our $\mathscr{P} _{_{H_N}}$ as
\begin{equation}\label{eq_appendix_1}
\begin{split}
\int_\mathbb{P}\min\limits_{u\in H_N}||u-\mathfrak{u}^p||_{H}^2\,d\mu(p)=& \int_\mathbb{P} ||(I-\mathscr{P} _{_{H_N}}){\mathfrak{u}^p}||_{H}^2\,d\mu(p)\\
=& \int_\mathbb{P} \langle (I-\mathscr{P} _{_{H_N}}){\mathfrak{u}^p},{\mathfrak{u}^p}\rangle_H\,d\mu(p)\\
=& \text{Tr}\left((I-\mathscr{P} _{_{H_N}})\int_\mathbb{P} {\mathfrak{u}^p}\otimes {\mathfrak{u}^p}\,d\mu(p)\right).
\end{split}
\end{equation}
{Herein, $\text{Tr}$ represents the trace of the operator $(I-\mathscr{P} _{_{H_N}})A$, defined on a Hilbert space as
\begin{equation}
\text{Tr}((I-\mathscr{P} _{_{H_N}})A) = \sum _{i = 1}^\infty \langle(I-\mathscr{P} _{_{H_N}})A e_i,e_i\rangle_H,
\end{equation}
where $e_i$ is any orthonormal basis of $H$. This definition is independent of the chosen basis. In particular, if the operator admits an eigenbasis, the trace equals the sum of all eigenvalues (provided these are finite).} We remark that \eqref{eq_operatpr_A} defines a compact, symmetric, and positive semidefinite operator, as it can be approximated by a sum of symmetric, rank-one, positive semidefinite operators.
In particular, by classical spectral theory \cite{ciarlet2013linear}, we have that there exists a non-increasing sequence of eigenvalues $\lambda_k\to 0$ and an orthonormal eigenbasis $\phi_k$ such that 
\begin{equation}
\int_\mathbb{P} \mathfrak{u}^p\otimes \mathfrak{u}^p\,d\mu(p)=\sum\limits_{k=1}^\infty \lambda_k\phi_k\otimes \phi_k. 
\end{equation}
Furthermore, it is a trace-class operator. This can be seen by taking an arbitrary orthonormal basis $e_i$ of $H$, and then we have that 
\begin{equation}
\begin{split}
\displaystyle \sum\limits_{k=1}^\infty \lambda_k = &\sum\limits_{i=1}^\infty \left\langle e_i,\int_\mathbb{P} \mathfrak{u}^p\otimes \mathfrak{u}^p\,d\mu(p)e_i\right\rangle_H\\
=&\int_\mathbb{P}\sum\limits_{i=1}^\infty \left\langle e_i, (\mathfrak{u}^p\otimes \mathfrak{u}^p) e_i\right\rangle_H\,d\mu(p)\\
=&\int_\mathbb{P}\sum\limits_{i=1}^\infty \langle \mathfrak{u}^p,e_i\rangle^2_H\,d\mu(p)\\
=&\int_\mathbb{P}||\mathfrak{u}^p||^2_H\,d\mu(p)<+\infty.
\end{split}
\end{equation}
Therefore, \eqref{eq_appendix_1} can be written as 
\begin{equation}
\begin{split}
\int_\mathbb{P} \min\limits_{u\in H_N}||u-\mathfrak{u}^p||_H^2\,d\mu(p)= & \text{Tr}\left((I-\mathscr{P} _{_{H_N}})\sum\limits_{k=1}^\infty \lambda_k\phi_k\otimes \phi_k\right)\\
=& \sum\limits_{k=1}^\infty \lambda_k ||(I-\mathscr{P} _{_{H_N}})\phi_k||_H^2.
\end{split}
\end{equation}
As $\lambda_k$ is summable, this means that for $N$ sufficiently large, $\displaystyle \sum\limits_{k=N+1}^\infty\lambda_k<\epsilon$, so by taking $H_N$ to be the span of $\{\phi_k\}_{k=1}^N$, we have that $(I-\mathscr{P} _{_{H_N}})\phi_k=\phi_k$ for $k\geq N$, and $(I-\mathscr{P} _{_{H_N}})\phi_k=0$ otherwise, so that
\begin{equation}
\begin{split}
\int_{\mathbb{P}} \min\limits_{u\in H_N}||u-\mathfrak{u}^p||_H^2 {d\mu(p)}= & \sum\limits_{k=1}^\infty \lambda_k ||(I-\mathscr{P} _{_{H_N}})\phi_k||_H^2\\
=& \sum\limits_{k=N+1}^\infty \lambda_k||\phi_k||_H^2\\
=& \sum\limits_{k=N+1}^\infty \lambda_k<\epsilon,
\end{split}
\end{equation}
from which we obtain the final result.
\end{proof}

\begin{theorem}\label{thm.UniversalApproxL2}
For every $\epsilon>0$, there exists a fully-connected feed-forward NN $\mathbf{u}^\alpha:\mathbb{R}^{N_0}\to\mathbb{R}^N$ such that if $H_N^\alpha=\text{span}\{u_n^\alpha\}_{n=1}^N$, then 
\begin{equation}
\int_\mathbb{P}||u^{p, \alpha}-\mathfrak{u}^p||_H^2\,d\mu(p)<\epsilon.
\end{equation}
\end{theorem}
\begin{proof}
First, we take $H_N$ to be a finite-dimensional subspace so that 
\begin{equation}
\int_\mathbb{P}\min\limits_{u\in H_N}||u-\mathfrak{u}^p||_H^2\,d\mu(p)<\frac{\epsilon}{4},
\end{equation}
which exists due to the uniform estimate in \Cref{prop.AppSubspaceL2}.
We denote 
\begin{equation}
C=\max\limits_{ p\in\mathbb{P}}||\mathfrak{u}^p||_H.
\end{equation}
As $\mathbb{P}$ is compact and $p\mapsto \mathfrak{u}^p$ is continuous, $C$ is necessarily finite. We take $\mathbf{u}^\alpha$ to be a fully-connected feed-forward NN such that $\displaystyle \left(\sum\limits_{n=1}^N ||u^\alpha_n-u_n||_{H}\right)^2<\frac{\epsilon}{4C^2}$, where $u_n$ is an orthonormal basis of $H_N$, and whose existance is guaranteed due to \Cref{prop.AppNNSubspace}. Given $p\in \mathbb{P}$, we define the coefficients $c_n(p) = \langle \mathfrak{u}^p,u_n\rangle_H$ and $ \tilde{u}^p = \displaystyle \sum\limits_{n=1}^Nc_nu_n^\alpha$. We observe that $c_n$ are uniformly bounded, as
\begin{equation}
|c_n(\xi)|\leq ||\mathfrak{u}^p||_H||u_n||_H=||\mathfrak{u}^p||_H\leq C. 
\end{equation}
Therefore, we have that 
\begin{equation}
\begin{split}
\int_\mathbb{P}||u^{p, \alpha}-\mathfrak{u}^p||_H^2\,d\mu(p)\leq & \int_\mathbb{P} ||\tilde{u}^p -\mathfrak{u}^p||_H^2\,d\mu(p)\\
{=} & \int_\mathbb{P}  \left|\left|\sum\limits_{n=1}^N c_n(u^\alpha_n-u_n)-(I-\mathscr{P} _{_{H_N}})\mathfrak{u}^p\right|\right|_{H}^2\,d\mu(p)\\
\leq & 2\int_\mathbb{P}\left|\left|\sum\limits_{n=1}^N c_n(u^\alpha_n-u_n)\right|\right|_{H}^2+\left|\left|(I-\mathscr{P} _{_{H_N}})\mathfrak{u}^p\right|\right|_{H}^2\,d\mu(p)\\
\leq &2\int_\mathbb{P}\left(\sum\limits_{n=1}^N\left|\left| c_n(u^\alpha_n-u_n)\right|\right|_{H}\right)^2\,d\mu(p)+2\int_\mathbb{P}\left|\left|(I-\mathscr{P} _{_{H_N}})\mathfrak{u}^p\right|\right|_{H}^2\,d\mu(p)\\
\leq & 2C^2 \int_\mathbb{P}\left(\sum\limits_{n=1}^N\left|\left|(u^\alpha_n-u_n)\right|\right|_{H}\right)^2\,d\mu(p)+2\int_\mathbb{P}\left|\left|(I-\mathscr{P} _{_{H_N}})\mathfrak{u}^p\right|\right|_{H}^2\,d\mu(p)\\
\leq & 2C^2\frac{\epsilon}{4C^2}+2\frac{\epsilon}{4}=\epsilon
\end{split}
\end{equation}

\end{proof}

\begin{corollary}
For every $\epsilon>0$, there exists a fully-connected feed-forward NN $\mathbf{u}^\alpha:\mathbb{R}^{N_0}\to\mathbb{R}^N$ such that if $H_N^\alpha=\text{span}\{u_n\}_{n=1}^N$, then 
\begin{equation}
\mathcal{L}_\mu(\alpha) <\epsilon
\end{equation}

\end{corollary}
\begin{proof}
This is a direct consequence of \Cref{thm.UniversalApproxL2} and the uniform estimates \eqref{eq_Continuity condition} and \eqref{eq_Inf-sup stability condition}.
\end{proof}

\begin{remark}
We defined a compact, trace-class, symmetric, and positive semi-definite operator \eqref{eq_operatpr_A} in the proof of \Cref{prop.AppSubspaceL2}. It was shown that the minimizer of
\begin{equation}\label{eq.min.V.norm}
\int_{\mathbb{P}} \min\limits_{u \in H_N}||u-\mathfrak{u}^p||_{H}^2\,d\mu(p)
\end{equation}
over all subspaces $H_N \subset H$ of a given dimension $N$ corresponds to the span of the first $N$ eigenvectors of operator $A$, with the error being 
\begin{equation}
\sum\limits_{k=N+1}^\infty \lambda_k. 
\end{equation}
Thus, we can understand the subspace $H_N$ of a given dimension that minimizes the integral \eqref{eq.min.V.norm} to correspond to the space spanned by the first $N$ eigenvectors of an operator that acts like an averaged projection onto the space of all possible solutions $\mathfrak{u}^p$.  {This approach is conceptually similar to the intuition behind the PGD method, which also relies on a projection operator defined with respect to the $L^2$-norm to construct an optimal subspace for approximating the solution space.} Essentially, this result, in the case that $\vartheta=\gamma$, tells us that we are searching for an optimal $H^k$-based projection operator onto the solution space. In the case where $\vartheta>\gamma$, this interpretation is less clear, however. 
\end{remark}

\bibliographystyle{abbrv}
\bibliography{references}

\begin{thebibliography}{10}

\bibitem{alvarez2013inversion}
J.~Alvarez-Aramberri, D.~Pardo, and H.~Barucq.
\newblock Inversion of magnetotelluric measurements using multigoal oriented
  hp-adaptivity.
\newblock {\em Procedia Computer Science}, 18:1564--1573, 2013.

\bibitem{azeez2001proper}
M.~Azeez and A.~Vakakis.
\newblock Proper orthogonal decomposition ({POD}) of a class of vibroimpact
  oscillations.
\newblock {\em Journal of Sound and vibration}, 240(5):859--889, 2001.

\bibitem{babuvska1971error}
I.~Babu{\v{s}}ka.
\newblock Error-bounds for finite element method.
\newblock {\em Numerische Mathematik}, 16(4):322--333, 1971.

\bibitem{bachmayr2024variationally}
M.~Bachmayr, W.~Dahmen, and M.~Oster.
\newblock Variationally correct neural residual regression for parametric
  {PDE}s: On the viability of controlled accuracy.
\newblock {\em arXiv preprint arXiv:2405.20065}, 2024.

\bibitem{software_LS-Net4ParametricPDEs}
S.~Baharlouei and C.~Uriarte.
\newblock {LS-Net4ParametricPDEs}.
\newblock \url{https://github.com/Mathmode/LS-Net4ParametricPDEs}, 2024.

\bibitem{berrone2023enforcing}
S.~Berrone, C.~Canuto, M.~Pintore, and N.~Sukumar.
\newblock Enforcing {D}irichlet boundary conditions in physics-informed neural
  networks and variational physics-informed neural networks.
\newblock {\em Heliyon}, 9(8), 2023.

\bibitem{bhattacharya2021model}
K.~Bhattacharya, B.~Hosseini, N.~B. Kovachki, and A.~M. Stuart.
\newblock Model reduction and neural networks for parametric {PDE}s.
\newblock {\em The SMAI journal of computational mathematics}, 7:121--157,
  2021.

\bibitem{bonilla2008inverse}
L.~L. Bonilla, A.~Carpio, O.~Dorn, M.~Moscoso, F.~Natterer, G.~Papanicolaou,
  M.~Rap{\'u}n, and A.~Teta.
\newblock {\em Inverse Problems and Imaging}.
\newblock Springer, 2008.

\bibitem{brevis2024learning}
I.~Brevis, I.~Muga, D.~Pardo, O.~Rodriguez, and K.~G. Van Der~Zee.
\newblock Learning quantities of interest from parametric {PDE}s: An efficient
  neural-weighted minimal residual approach.
\newblock {\em Computers \& Mathematics with Applications}, 164:139--149, 2024.

\bibitem{chen2014inf}
L.~Chen.
\newblock Inf-sup conditions.
\newblock {\em University of California at Irvine, Irvine, CA, Technical
  report}, 2014.

\bibitem{chen1995approximation}
T.~Chen and H.~Chen.
\newblock Approximation capability to functions of several variables, nonlinear
  functionals, and operators by radial basis function neural networks.
\newblock {\em IEEE Transactions on Neural Networks}, 6(4):904--910, 1995.

\bibitem{chen1995universal}
T.~Chen and H.~Chen.
\newblock Universal approximation to nonlinear operators by neural networks
  with arbitrary activation functions and its application to dynamical systems.
\newblock {\em IEEE transactions on neural networks}, 6(4):911--917, 1995.

\bibitem{chinesta2011overview}
F.~Chinesta, A.~Ammar, A.~Leygue, and R.~Keunings.
\newblock An overview of the proper generalized decomposition with applications
  in computational rheology.
\newblock {\em Journal of Non-Newtonian Fluid Mechanics}, 166(11):578--592,
  2011.

\bibitem{chinesta2011short}
F.~Chinesta, P.~Ladeveze, and E.~Cueto.
\newblock A short review on model order reduction based on proper generalized
  decomposition.
\newblock {\em Archives of Computational Methods in Engineering},
  18(4):395--404, 2011.

\bibitem{ciarlet2013linear}
P.~G. Ciarlet.
\newblock {\em Linear and nonlinear functional analysis with applications}.
\newblock SIAM, 2013.

\bibitem{dal2020data}
N.~Dal~Santo, S.~Deparis, and L.~Pegolotti.
\newblock Data driven approximation of parametrized {PDE}s by reduced basis and
  neural networks.
\newblock {\em Journal of Computational Physics}, 416:109550, 2020.

\bibitem{discacciati2024overlapping}
M.~Discacciati, B.~J. Evans, and M.~Giacomini.
\newblock An overlapping domain decomposition method for the solution of
  parametric elliptic problems via proper generalized decomposition.
\newblock {\em Computer Methods in Applied Mechanics and Engineering},
  418:116484, 2024.

\bibitem{garcia2017monitoring}
R.~Garc{\'\i}a-Blanco, D.~Borzacchiello, F.~Chinesta, and P.~Diez.
\newblock Monitoring a {PGD} solver for parametric power flow problems with
  goal-oriented error assessment.
\newblock {\em International Journal for Numerical Methods in Engineering},
  111(6):529--552, 2017.

\bibitem{geist2021numerical}
M.~Geist, P.~Petersen, M.~Raslan, R.~Schneider, and G.~Kutyniok.
\newblock Numerical solution of the parametric diffusion equation by deep
  neural networks.
\newblock {\em Journal of Scientific Computing}, 88(1):22, 2021.

\bibitem{goswami2022physics}
S.~Goswami, M.~Yin, Y.~Yu, and G.~E. Karniadakis.
\newblock A physics-informed variational deeponet for predicting crack path in
  quasi-brittle materials.
\newblock {\em Computer Methods in Applied Mechanics and Engineering},
  391:114587, 2022.

\bibitem{haasdonk2008reduced}
B.~Haasdonk and M.~Ohlberger.
\newblock Reduced basis method for finite volume approximations of parametrized
  linear evolution equations.
\newblock {\em ESAIM: Mathematical Modelling and Numerical Analysis},
  42(2):277--302, 2008.

\bibitem{han2018solving}
J.~Han, A.~Jentzen, and W.~E.
\newblock Solving high-dimensional partial differential equations using deep
  learning.
\newblock {\em Proceedings of the National Academy of Sciences},
  115(34):8505--8510, 2018.

\bibitem{hille1996functional}
E.~Hille and R.~S. Phillips.
\newblock {\em Functional analysis and semi-groups}, volume~31.
\newblock American Mathematical Soc., 1996.

\bibitem{hornik1990universal}
K.~Hornik, M.~Stinchcombe, and H.~White.
\newblock Universal approximation of an unknown mapping and its derivatives
  using multilayer feedforward networks.
\newblock {\em Neural networks}, 3(5):551--560, 1990.

\bibitem{khara2024neufenet}
B.~Khara, A.~Balu, A.~Joshi, S.~Sarkar, C.~Hegde, A.~Krishnamurthy, and
  B.~Ganapathysubramanian.
\newblock Neufenet: Neural finite element solutions with theoretical bounds for
  parametric {PDE}s.
\newblock {\em Engineering with Computers}, pages 1--23, 2024.

\bibitem{kharazmi2019variational}
E.~Kharazmi, Z.~Zhang, and G.~E. Karniadakis.
\newblock Variational physics-informed neural networks for solving partial
  differential equations.
\newblock {\em arXiv preprint arXiv:1912.00873}, 2019.

\bibitem{khoo2021solving}
Y.~Khoo, J.~Lu, and L.~Ying.
\newblock Solving parametric {PDE} problems with artificial neural networks.
\newblock {\em European Journal of Applied Mathematics}, 32(3):421--435, 2021.

\bibitem{khoo2019switchnet}
Y.~Khoo and L.~Ying.
\newblock Switchnet: a neural network model for forward and inverse scattering
  problems.
\newblock {\em SIAM Journal on Scientific Computing}, 41(5):A3182--A3201, 2019.

\bibitem{khoromskij2011tensor}
B.~N. Khoromskij and C.~Schwab.
\newblock Tensor-structured galerkin approximation of parametric and stochastic
  elliptic {PDE}s.
\newblock {\em SIAM journal on scientific computing}, 33(1):364--385, 2011.

\bibitem{kovachki2023neural}
N.~Kovachki, Z.~Li, B.~Liu, K.~Azizzadenesheli, K.~Bhattacharya, A.~Stuart, and
  A.~Anandkumar.
\newblock Neural operator: Learning maps between function spaces with
  applications to {PDE}s.
\newblock {\em Journal of Machine Learning Research}, 24(89):1--97, 2023.

\bibitem{kutyniok2022theoretical}
G.~Kutyniok, P.~Petersen, M.~Raslan, and R.~Schneider.
\newblock A theoretical analysis of deep neural networks and parametric {PDE}s.
\newblock {\em Constructive Approximation}, 55(1):73--125, 2022.

\bibitem{li2020fourier}
Z.~Li, N.~Kovachki, K.~Azizzadenesheli, B.~Liu, K.~Bhattacharya, A.~Stuart, and
  A.~Anandkumar.
\newblock Fourier neural operator for parametric partial differential
  equations.
\newblock {\em arXiv preprint arXiv:2010.08895}, 2020.

\bibitem{li2020neural}
Z.~Li, N.~Kovachki, K.~Azizzadenesheli, B.~Liu, K.~Bhattacharya, A.~Stuart, and
  A.~Anandkumar.
\newblock Neural operator: Graph kernel network for partial differential
  equations.
\newblock {\em arXiv preprint arXiv:2003.03485}, 2020.

\bibitem{lu2019deeponet}
L.~Lu, P.~Jin, and G.~E. Karniadakis.
\newblock Deeponet: Learning nonlinear operators for identifying differential
  equations based on the universal approximation theorem of operators.
\newblock {\em arXiv preprint arXiv:1910.03193}, 2019.

\bibitem{lu2021learning}
L.~Lu, P.~Jin, G.~Pang, Z.~Zhang, and G.~E. Karniadakis.
\newblock Learning nonlinear operators via deeponet based on the universal
  approximation theorem of operators.
\newblock {\em Nature machine intelligence}, 3(3):218--229, 2021.

\bibitem{lu2022comprehensive}
L.~Lu, X.~Meng, S.~Cai, Z.~Mao, S.~Goswami, Z.~Zhang, and G.~E. Karniadakis.
\newblock A comprehensive and fair comparison of two neural operators (with
  practical extensions) based on fair data.
\newblock {\em Computer Methods in Applied Mechanics and Engineering},
  393:114778, 2022.

\bibitem{nelsen2021random}
N.~H. Nelsen and A.~M. Stuart.
\newblock The random feature model for input-output maps between banach spaces.
\newblock {\em SIAM Journal on Scientific Computing}, 43(5):A3212--A3243, 2021.

\bibitem{nouy2010priori}
A.~Nouy.
\newblock A priori model reduction through proper generalized decomposition for
  solving time-dependent partial differential equations.
\newblock {\em Computer Methods in Applied Mechanics and Engineering},
  199(23-24):1603--1626, 2010.

\bibitem{pinkus1999approximation}
A.~Pinkus.
\newblock Approximation theory of the mlp model in neural networks.
\newblock {\em Acta numerica}, 8:143--195, 1999.

\bibitem{raissi2017physics}
M.~Raissi, P.~Perdikaris, and G.~E. Karniadakis.
\newblock Physics informed deep learning (part i): Data-driven solutions of
  nonlinear partial differential equations.
\newblock {\em arXiv preprint arXiv:1711.10561}, 2017.

\bibitem{raissi2019physics}
M.~Raissi, P.~Perdikaris, and G.~E. Karniadakis.
\newblock Physics-informed neural networks: A deep learning framework for
  solving forward and inverse problems involving nonlinear partial differential
  equations.
\newblock {\em Journal of Computational physics}, 378:686--707, 2019.

\bibitem{rojas2024robust}
S.~Rojas, P.~Maczuga, J.~Mu{\~n}oz-Matute, D.~Pardo, and M.~Paszy{\'n}ski.
\newblock Robust variational physics-informed neural networks.
\newblock {\em Computer Methods in Applied Mechanics and Engineering},
  425:116904, 2024.

\bibitem{rozza2008reduced}
G.~Rozza, D.~B.~P. Huynh, and A.~T. Patera.
\newblock Reduced basis approximation and a posteriori error estimation for
  affinely parametrized elliptic coercive partial differential equations:
  application to transport and continuum mechanics.
\newblock {\em Archives of Computational Methods in Engineering},
  15(3):229--275, 2008.

\bibitem{shukla2024deep}
K.~Shukla, V.~Oommen, A.~Peyvan, M.~Penwarden, N.~Plewacki, L.~Bravo,
  A.~Ghoshal, R.~M. Kirby, and G.~E. Karniadakis.
\newblock Deep neural operators as accurate surrogates for shape optimization.
\newblock {\em Engineering Applications of Artificial Intelligence},
  129:107615, 2024.

\bibitem{sibileau2018explicit}
A.~Sibileau, A.~Garc{\'\i}a-Gonz{\'a}lez, F.~Auricchio, S.~Morganti, and
  P.~D{\'\i}ez.
\newblock Explicit parametric solutions of lattice structures with proper
  generalized decomposition ({PGD}) applications to the design of 3d-printed
  architectured materials.
\newblock {\em Computational Mechanics}, 62:871--891, 2018.

\bibitem{taylor2023deep}
J.~M. Taylor, D.~Pardo, and I.~Muga.
\newblock A deep {F}ourier residual method for solving {PDE}s using neural
  networks.
\newblock {\em Computer Methods in Applied Mechanics and Engineering},
  405:115850, 2023.

\bibitem{taylor2024regularity}
J.~M. Taylor, D.~Pardo, and J.~Mu{\~n}oz-Matute.
\newblock Regularity-conforming neural networks ({ReCoNNs}) for solving partial
  differential equations.
\newblock {\em arXiv preprint arXiv:2405.14110}, 2024.

\bibitem{uriarte2024optimizing}
C.~Uriarte, M.~Bastidas, D.~Pardo, J.~M. Taylor, and S.~Rojas.
\newblock Optimizing variational physics-informed neural networks using least
  squares.
\newblock {\em arXiv preprint arXiv:2407.20417}, 2024.

\bibitem{uriarte2023deep}
C.~Uriarte, D.~Pardo, I.~Muga, and J.~Mu{\~n}oz-Matute.
\newblock A {D}eep {D}ouble {R}itz {M}ethod ({D2RM}) for solving {P}artial
  {D}ifferential {E}quations using {N}eural {N}etworks.
\newblock {\em Computer Methods in Applied Mechanics and Engineering},
  405:115892, 2023.

\bibitem{uriarte2022finite}
C.~Uriarte, D.~Pardo, and {\'A}.~J. Omella.
\newblock A finite element based deep learning solver for parametric {PDE}s.
\newblock {\em Computer Methods in Applied Mechanics and Engineering},
  391:114562, 2022.

\bibitem{wang2023long}
S.~Wang and P.~Perdikaris.
\newblock Long-time integration of parametric evolution equations with
  physics-informed {DeepONets}.
\newblock {\em Journal of Computational Physics}, 475:111855, 2023.

\bibitem{wang2021learning}
S.~Wang, H.~Wang, and P.~Perdikaris.
\newblock Learning the solution operator of parametric partial differential
  equations with physics-informed {DeepONets}.
\newblock {\em Science advances}, 7(40):eabi8605, 2021.

\bibitem{yu2018deep}
B.~Yu et~al.
\newblock The deep {Ritz} method: a deep learning-based numerical algorithm for
  solving variational problems.
\newblock {\em Communications in Mathematics and Statistics}, 6(1):1--12, 2018.

\bibitem{zhdanov2002geophysical}
M.~S. Zhdanov.
\newblock {\em Geophysical inverse theory and regularization problems},
  volume~36.
\newblock Elsevier, 2002.

\end{thebibliography}
\end{document}